\documentclass[12pt]{article}
\usepackage{mathrsfs}
\usepackage{amsfonts}
\usepackage{amssymb}
\usepackage{amsthm}
\usepackage{amsmath}
\usepackage{latexsym}
\usepackage{color}

\allowdisplaybreaks[4]
\newtheorem{lemma}{\qquad Lemma}

\newtheorem{remark}{\qquad Remark}

\newtheorem{thm}{\qquad Theorem}

\begin{document}

\title{\large\bf The uniform local asymptotics of the total net loss process in a new time-dependent bidimensional renewal model
\thanks{Research supported by the National Science Foundation of China
(No. 11071182 {\&} No. 71171177), the project of the
key research base of human and social
science (Statistics, Finance) for colleges in Zhejiang Province (Grant No. of
Academic Education of Zhejiang,~2008-255).}}
\author {Tao Jiang$^{1)}$\qquad Yuebao Wang$^{2)}$ \thanks{Corresponding author. \qquad Telephone: +86 512 67422726. \qquad Fax: +86 512 65112637. \quad E-mail: ybwang@suda.edu.cn (Y. Wang)}
\qquad Hui Xu$^{2)}$ \qquad\\ {\small\it 1). School of Statistics and Mathematics, Zhejiang Gongshang University, P. R. China, 310018 }\\ {\small\it 2). School of Mathematical Sciences, Soochow University, Suzhou, P. R. China, 215006} \\
}
\date{}

\maketitle

\begin{abstract}
In this paper, we consider a bidimensional renewal risk model with constant force of interest, in which the claim size vector with certain local subexponential marginal distribution and its inter-arrival time are subject to a new time-dependence structure. We obtain the uniform local asymptotics of the total net loss process in the model. Moreover, some specific examples of the joint  distribution satisfying the conditions of the dependence structure are given. Finally, in order to illustrate a condition of the above result,
a local subexponential distribution is find for the first time that, its local distribution is not almost decreased.\\
\medskip

{\it Keywords:} uniform local asymptotics; total net loss process; time dependence; bidimensional renewal model; local subexponential distribution; almost decreasing\\

{\it 2000 Mathematics Subject Classification:} Primary 62H05; 62E20; 62P05
\end{abstract}

\section{Introduction}
It is well known that, since Kl\"{u}ppelberg and Stadtimuller (1998) began to study the unidimensional renewal risk models with constant force of interest, there are many related researches; see, for example, Kalashnikov and Konstantinides (2000), Konstantinides et al. (2002), Tang (2005, 2007), and Hao and Tang (2008). In these works, the claim sizes
$X_k,k\in\mathbb{N}$ are usually assumed to be a sequence of independent and
identically distributed (i.i.d.) random variables (r.v.$^,$s) with
generic r.v. $X$, and their inter-arrival times $\theta_k,k\in\mathbb{N}$ are a sequence of i.i.d. r.v.$^,$s
with generic r.v. $\theta$. Furthermore, unidimensional renewal risk models have also been investigated with certain dependence structure among either the claim sizes or the inter-arrival times; see Chen and Ng (2007), Liu et al. (2012), Wang et al. (2013), etc. In these papers, the claim sizes and their inter-arrival times are assumed mutually independent,
and the imposed dependence structure among either the claim sizes or the
inter-arrival times yields have no impact on the asymptotic behavior of certain research object.

The research on a risk model with a certain dependence structure between $X$
and $\theta$ can be found in Albrecher and Teugels (2006), Badescu et al. (2009),
Asimit and Badescu (2010), and Li et al. (2010). These papers show that the object of study indeed depends on the proposed dependence structure between the claim size and inter-arrival times.

The literatures on multidimensional risk models focused on a continuous-time setting;  see Chan et al. (2003), Li et al. (2007), Chen et al. (2011),  Chen et al. (2013b), Yang and Li (2014), Jiang et al. (2015), and so on. A multivariate risk model in a discrete-time framework was studied by Huang et al. (2014), etc.

So far, we have not found any results concerning local asymptotics for certain object of study in the multidimensional risk models. In the reality of insurance with heavy-tailed claim
sizes, the local probability of certain object of study is often a infinite small amount of the corresponding global probability. Therefore, the study of the local probability is of great importance both theoretically and pragmatically.

In this paper, we give the uniform asymptotic estimates for the local probabilities of the total net loss process in a new time-dependent bidimesional risk model for presentational convenience. The result can be trivially extended to the multidimensional case. To this end, in the following, we introduce the bidimesional renewal risk model, the local distribution classes, the dependent structure between the random vector of claim sizes and the inter-arrival time, and the main result of present paper, respectively.

\subsection{A bidimesional renewal risk model}

Assume that the insurance company has two classes of business. For $i=1,2$,
the claim sizes for the $i$-th class
$X_{k}^{(i)}, k\in\mathbb{N}$ are i.i.d. r.v.$^,$s
with common continuous distribution $F_i$ supported on $[0,\infty)$
that is $F_i(x)<1$ for all $x\ge 0$. We also assume that
the two claim sizes $X_{k}^{(1)}$ and $X_{k}^{(2)}$
occur at the same time for all
$k\in\mathbb{N}$, and their inter-arrival times
$\{\theta_k:k\in\mathbb{N}\}$ form another sequence of i.i.d. r.v.$^,$s
with common continuous distribution $G$ supported
on $[0,\infty)$.
The arrival times of successive claims are defined
by $\sigma_n=\sum_{k=1}^n\theta_k,n\in\mathbb{N}$, constituting an
ordinary renewal counting process
\begin{equation*}
N(t)=\sum_{n=1}^{\infty}\textbf{1}_{\{\sigma_n\leq t\}},~t\geq0.
\end{equation*}
Denote the renewal function by $\lambda(t)=EN(t)$ for all $0<t<\infty$ and $\lambda(0)=0$.
Define
$$\Lambda=\{t: \lambda(t)>0\}=\{t: P(\sigma_1\leq t)>0\}$$
for later use. With $\underline{t}=\inf\{t:P(\sigma_1\leq t)>0\} $, clearly, $\Lambda=[\underline{t},\infty]$ if
$P(\sigma_1=\underline{t})>0$; or $\Lambda=(\underline{t},\infty]$
if $P(\sigma_1=\underline{t})=0$.

Let $r\geq0$ be the constant force of interest.
The premiums accumulated up to time $t$ for the $i$th class, denoted by $C_i(t)$ with $C_i(0)=0$ and $C_i(t)<\infty$ almost surely (a.s.), $i=1,2$, follows two nonnegative and nondecreasing stochastic processes. Denote the vector by $\vec{C}(t)=(C_1(t),C_2(t))^{\top}$. Let $\vec{X}_k=(X_k^{(1)},X_k^{(2)})^{\top}$ be the $k$-th pair of claims, $k\ge1$, and $\vec{x}=(x_1,x_2)^{\top}$ be the initial surplus vector. The total net surplus up to $t$, denoted by $\vec{U}(\vec{x},t)=(U_{1}(x_1,t),U_{2}(x_2,t))^{\top}$, satisfies
\begin{equation}\label{p2}
\vec{U}(\vec{x},t)=\vec{x}e^{rt}+\int_{0-}^te^{r(t-y)}\vec{C}(dy)-\int_{0-}^te^{r(t-y)}\vec{S}(dy),
\end{equation}
where $\vec{S}(t)=\sum_{k=1}^{N(t)}\vec{X}_k,~t\geq0$. In the following, the $\{\vec{U}(\vec{x},t):t\in\Lambda\}$ is called the total net surplus process, and the $\{-\vec{U}(\vec{x},t):t\in\Lambda\}$ is called the total net loss process.
The vector of discounted aggregate claims is expressed as
\begin{eqnarray}\label{p3}
\vec{D}_r(t)&=&\int_{0-}^te^{-ry}\vec{S}(dy)=\sum_{k=1}^{\infty}\vec{X}_ke^{-r\sigma_k}1_{\{\sigma_k\leq t\}}\nonumber\\
&=&\sum_{n=1}^\infty\Big(\sum_{k=1}^n\vec{X}_ke^{-r\sigma_k}\Big)1_{\{N(t)=n\}}=\sum_{k=1}^{N(t)}\vec{X}_ke^{-r\sigma_k}.
\end{eqnarray}
Here and thereafter, for vectors $\vec{a}=(a_1,a_2)^{\top}$ and $\vec{b}=(b_1,b_2)^{\top}$, we write $\vec{a}\leq\vec{b}$ if
$a_1\leq b_1$ and $a_2\leq b_2$, write $\vec{a}<\vec{b}$ if $a_1<
b_1$ and $a_2< b_2$.

\subsection{Some local distribution classes}

For any constant $d\in(0,\infty]$ and some distribution $F$ supported on $[0,\infty)$, denote $F(x+\Delta_d)=F(x,x+d]$ when $d<\infty$, and $F(x+\Delta_d)=F(x,x+\infty)=\overline{F}(x)$ when $d=\infty$.

We say that a distribution $F$ belongs to the distribution class $\mathcal{L}_{loc}$, if for all $d\in(0,\infty],F(x+\Delta_d)>0$ eventually, and for any $t>0$ it holds uniformly for all $\mid s\mid\le t$ that
$$F(x+s+\Delta_d)\sim F(x+\Delta_d),$$
where $a(x)\sim b(x)$ for two positive functions $a(\cdot)$ and $b(\cdot)$ whenever
$\lim a(x)(b(x))^{-1}=1$, and all limit relationships are for $x\rightarrow\infty$,
unless otherwise stated. Clearly, the distribution $F$ in the class $\mathcal{L}_{loc}$ is heavy-tailed, that is $\int_{0}^{\infty}\exp\{\varepsilon y\}F(dy)=\infty$ holds for any $\varepsilon>0$

Further, if a distribution $F$ belongs to the class $\mathcal{L}_{loc}$ and
$$F^{*2}(x+\Delta_d)\sim 2F(x+\Delta_d),$$
then we say that the distribution $F$ belongs to the distribution class $\mathcal{S}_{loc}$, where $F^{*n}$ is the $n$-th convolution of $F$ with itself for $n\ge2$ and $F^{*1}=F$. See, for example, Borovkov and Borovkov \cite{BB2008}.

In aforementioned two conceptions, we replace ``for all $d$" with ``for some $d$", then we say that the distribution $F$ belongs to the local long-tailed distribution class $\mathcal{L}_{\Delta_d}$ and local subexponential distribution class $\mathcal{S}_{\Delta_d}$, respectively. See Asmussen et al. \cite{AFK2003}.

These local distribution classes play a crucial role in the research of the local asymptotics of some studied objects. On the research concerning independent r.v.'s, besides the aforementioned papers, the readers can refer to Wang et al. (2005), Shneer
(2006), Wang et al. (2007), Denisov and Shneer (2007), Denisov et al. (2008), Cui et al. (2009), Chen et al. (2009), Yu et al. (2010), Watanabe and Yamamuro (2010), Yang et al. (2010), Wang and Wang (2011), Lin (2012), Wang et al. (2016), and so on. However, the research related to certain dependent r.v.'s is very rare.

\subsection{A new local time-dependent structure}

In this paper, we assume that $\{(X^{(1)}_k,X^{(2)}_k,\theta_k): k\in\mathbb{N}\}$
is a sequence of i.i.d. random vectors with generic vector $(X^{(1)},X^{(2)},\theta)$. Further, based on the idea of Asimit and Badescu (2010) for the global joint distribution, we construct a new dependence structure of $(X^{(1)},X^{(2)},\theta)$ satisfying the following conditions. Here, all the related functions are positive and measurable, and all limit relationships are for $(x_1,x_2)\to(\infty,\infty)$, unless otherwise stated. In addition, we denote $\Delta^{(i)}=(0,d_i]$ for $0<d_i<\infty,i=1,2$ and $T$ is any fixed positive number in $\Lambda$.

\vspace{0.2cm}
\ \ \ {\bf Condition 1.} For $i=1,2$, there exists a function $h_i(\cdot)$ such that
\begin{eqnarray}\label{202}
P(X^{(i)}\in x_i+\Delta^{(i)}\mid\theta=s)\sim F_i(x_i+\Delta^{(i)})h_i(s),
\end{eqnarray}
uniformly for all $s\ge0$ and $d_i>0$ satisfying $F_i(x_i+\Delta^{(i)})>0$; and
\begin{eqnarray}\label{203}
0<b_*&=&b_*(T)=\min\Big\{\inf_{s\in[0,T]}h_i(s):i=1,2\Big\}\nonumber\\
&\le&\max\Big\{\sup_{s\in[0,T]}h_i(s):\ i=1,2\Big\}=b^*(T)=b^*<\infty.
\end{eqnarray}

\ \ \ {\bf Condition 2.} There
exists a function $g(\cdot)$ such that
\begin{eqnarray}\label{204}
P(\vec{X}\in\vec{x}+\vec{\Delta} \mid\theta=s)\sim F_1(x_1+\Delta^{(1)})F_2(x_2+\Delta^{(2)})g(s),
\end{eqnarray}
uniformly for all $s\ge0$ and $d_i>0$ satisfying $F_i(x_i+\Delta^{(i)})>0$; and
\begin{eqnarray}\label{205}
0<d_*=d_*(T)=\inf_{s\in[0,T]}g(s)\le\sup_{s\in[0,T]}g(s)=d^*(T)=d^*<\infty.
\end{eqnarray}

\ \ \ {\bf Condition 3.} For $1\le i\neq j\le2$, there exists a binary function $g_{ij}(\cdot,\cdot)$ such that
\begin{eqnarray}\label{207}
P(X^{(i)}\in x_i+\Delta^{(i)}\mid X^{(j)}=z,\theta=s)\sim F_i(x_i+\Delta^{(i)})g_{ij}(z,s)
\end{eqnarray}
uniformly for all $z,s\ge0$ and $d_i>0$ satisfying $F_i(x_i+\Delta^{(i)})>0$; and
\begin{eqnarray}\label{208}
0<a_*&=&a_*(T)=\min\{\inf_{z\ge0,\ s\in[0,T]}g_{ij}(z,s):i=1,2,i\neq j\}\nonumber\\
&\le&\max\{\sup_{z\ge0,\ s\in[0,T]}g_{ij}(z,s):i=1,2,i\neq j\}=a^*(T)=a^*<\infty.
\end{eqnarray}

In the above,  when $s$ is not a possible value of $\theta$, the conditional probabilities should be understood as unconditional ones so that  $h_i(s)=g(s)=g_{ij}(z,s)=1$ for $1\le i\neq j\le2$ and $z\ge0$. In Section 4, some concrete copulas of $(X^{(1)},X^{(2)},\theta)$ satisfying the Conditions 1-3 will be presented. They include the well-known Sarmanov joint distribution (and hence the Farlie-Gumbel-Morgenstern joint distribution), Frank joint distribution, and other distributions. Interestingly, we embed a two-dimensional  product copula into a two-dimensional Frank copula, and we can get a such three-dimensional joint distribution.

The following condition is independent of the dependence structure of the model.

\ \ \ {\bf Condition 4.} For $i=1,2$, and for all $x_i,y_i,d_i>0$ satisfying $F_i(x_i+\Delta^{(i)})F_i(y_i+\Delta^{(i)})>0$, there is a constant $C_4=C_4(F_1,F_2)\ge0$ such that
\begin{eqnarray*}\label{2012}
1+C_{4}=\sup_{0\le x_i\le y_i<\infty,i=1,2}F_i(y_i+\Delta^{(i)})\big(F_i(x_i+\Delta^{(i)})\big)^{-1}<\infty.
\end{eqnarray*}

\begin{remark}
For $i=1,2$, according to Condition 1 and Condition 4, there is a constant $C_1=C_1(F_1,F_2)\ge0$ such that,
\begin{eqnarray}\label{201}
P(X^{(i)}\in x_i+\Delta^{(i)}\mid\theta=s)\le (1+C_1)F_i(x_i+\Delta^{(i)}),
\end{eqnarray}
for all $s,x_i\ge0$ and $d_i>0$ satisfying $F_i(x_i+\Delta^{(i)})>0$.

Similarly, according to Condition 2 and Condition 4, there is a constant $C_2=C_2(F_1,F_2)$ $\ge0$ such that,
\begin{eqnarray}\label{209}
P(\vec{X}\in\vec{x}+\vec{\Delta} \mid\theta=s)\le (1+C_2)F_1(x_1+\Delta^{(1)})F_2(x_2+\Delta^{(2)}),
\end{eqnarray}
for all $s, x_i\ge0$ and all $d_i>0$ satisfying $F_i(x_i+\Delta^{(i)})>0,i=1,2$.
Further, without loss of generality, we also assume that, for all $x_1,x_2,s\ge0$,
\begin{eqnarray}\label{2010}
P(X^{(1)}\in dx_1,X^{(2)}\in dx_2\mid\theta=s)\le (1+C_2)F_1(dx_1)F_2(dx_2)
\end{eqnarray}
and
\begin{eqnarray}\label{2011}
P(X^{(1)}\in dx_1,X^{(2)}\in dx_2,\theta\in ds)\le (1+C_2)F_1(dx_1)F_2(dx_2)G(ds).
\end{eqnarray}

Finally, according to Condition 3 and Condition 4, for $1\le i\neq j\le2$, there is a constant $C_3=C_3(F_1,F_2)\ge0$ such that
\begin{eqnarray}\label{206}
P(X^{(i)}\in x_i+\Delta^{(i)}\mid X^{(j)}=z,\theta=s)\le (1+C_3)F_i(x_i+\Delta^{(i)}).
\end{eqnarray}
for all $s,z, x_i\ge0$ and all $d_i>0$ satisfying $F_i(x_i+\Delta^{(i)})>0$.
\end{remark}

\begin{remark}
The Condition 4 is slightly stronger than the following condition:
\begin{eqnarray*}
\sup_{x_0\le x_i\le y_i<\infty,i=1,2}F_i(y_i+\Delta^{(i)})\big(F_i(x_i+\Delta^{(i)})\big)^{-1}<\infty
\end{eqnarray*}
for some positive constant $x_0=x_0(F_1,F_2)$. The condition in unidimensional case is used by Lemma 2.3 and Corollary 2.1 of Denisov et al (2008), Proposition 6.1 of Wang and Wang (2011), and so on. In the terminology of Bingham
et al. (1987), the condition is called that the local distribution of $F_i$ is almost decreased, or the distribution is locally almost decreased, for $i=1,2$.
And for many common distributions in $\mathcal{S}_{loc}$, their local distributions are almost decreased with $C_{4}=0$. However, there is a question that, are all local distributions in $\mathcal{S}_{loc}$ almost decreased? For the answer, we have not found any counter examples or positive proof in the literature. Therefore, in Section 5, we construct a distribution belonging to the class $\mathcal{S}_{loc}$, the local distribution of which is not almost decreased.
\end{remark}



\subsection{Main result}

Firstly, we denote the two-dimensional joint distributions of random vectors $\vec{D}_r(t)$ and $-\vec{U}(\vec{x},t)$ by $F_{\vec{D}_r(t)}$ and $F_{-\vec{U}(\vec{x},t)}$ for $t\in\Lambda$, respectively.
\begin{thm} \label{thm201} Consider above the bidimensional renewal risk model satisfying Conditions 1-4. Suppose that $F_i\in\mathcal{S}_{loc}$, $i=1,2$ and $Ee^{\beta N(T)}<\infty$ for some $\beta>C_2$ in (\ref{209}).
Then it holds uniformly for all $t\in \Lambda\cap[0,T]$ that
\begin{eqnarray}\label{th211}
&&F_{\vec{D}_r(t)}(\vec{x}+\vec{\Delta})\sim\int_{0-}^{t}\int_{0-}^{t-u}
\Big(P\big(X^{(1)}e^{-r(u+v)}\in x_1+\Delta^{(1)}\big)P\big(X^{(2)}e^{-rv}\in x_2+\Delta^{(2)}\big)\nonumber\\
&&\ \ \ \ \ \ \ \ \ +P\big(X^{(1)}e^{-ru}\in x_1+\Delta^{(1)}\big)P\big(X^{(2)}e^{-r(u+v)}\in x_2+\Delta^{(2)}\big)\Big)
\tilde{\lambda}_2(dv)\tilde{\lambda}_1(du)\nonumber\\
&&\ \ \ \ \ \ \ \ \ \ \ \ \ \ \ \ \ +\int_{0-}^{t}
\prod_{i=1}^2 P\big(X^{(i)}_ke^{-ru}\in x_i+\Delta^{(i)})\big)
\tilde{\tilde{\lambda}}(du),
\end{eqnarray}
where
$$\tilde{\lambda}_i(t)=\int_{0-}^{t}(1+\lambda(t-u))h_i(u)G(du),\ \ i=1,2,$$
and
$$\tilde{\tilde{\lambda}}(t)=\int_{0-}^{t}(1+\lambda(t-u))g(u)G(du).$$

Furthermore, assume that the processes $\{C_i(t):0\leq t<\infty, i=1,2\}$ are
independent of $\{X_{k}^{(i)}:k\in\mathbb{N},i=1,2\}$ and $\{N(t):t\geq
0\}$. Then the local probability of the total net loss process up to $t$
\begin{eqnarray}\label{th212}
P\big(-\vec{U}(\vec{x},t)\in(\vec{0},\vec{d}\ ]\big)=P\big(\vec{U}(\vec{x},t)\in[-\vec{d},\vec{0}\ )\big)\sim F_{\vec{D}_r(t)}(\vec{x}+e^{-rt}\vec{\Delta})
\end{eqnarray}
holds uniformly for all $t\in \Lambda\cap[0,T]$.
\end{thm}


\begin{remark}
It is well known that, for any fixed
$T\in\Lambda$, there is a constant $\beta_0=\beta_0(T)>0$ such that
$Ee^{\beta_0 N(T)}<\infty$. Thus, if such $\beta_0>C_2$, then the condition that
$Ee^{\beta N(T)}<\infty$ holds automatically for any $\beta\in(C_2,\beta_0)$. Particularly, when the dependency of $(X^{(1)},X^{(2)},\theta)$ is governed by a tri-dimensional Farlie-Gumbel-Morgenstern copula with $-1<\gamma_{12}\le0$, see Copula 4.1 below, we shall find that $(X^{(1)},X^{(2)},\theta)$ satisfies Condition 1-4
with $C_2=0$. Thus, we have $Ee^{\beta N(T)}<\infty$ for any $\beta\in(0, \beta_0)$.
\end{remark}

\begin{remark}
As what Li et al.(2010) discussed, for $i=1,2$, we
can introduce a r.v. $\theta^*_i$ with a proper distribution given by
\begin{equation}\label{rem211}
P(\theta^{*}_i\in dt)=\frac{h_i(t)}{Eh_i(\theta)\textbf{\emph{1}}_{\{\theta\leq
T\}}}G( dt),~~t\in[0,T].
\end{equation}
Let $\{\theta^*_{i,k}: k\in\mathbb{N}\}$ be
a sequence of i.i.d. r.v.'s with a distribution law specified in \eqref{rem211},
$i=1,2$. Then the inter-arrival times $\theta^*_{i,1},\theta_k, k\ge
2$, constitutes a delayed renewal counting process $\{N^*_i(t):
t\geq0\}$ with a corresponding renewal function
$\lambda^*_i(t)=EN^*_i(t),t\geq0,\ i=1,2$. It is easy to see that
\begin{equation}\label{rem212}
\tilde{\lambda}_i(dt)=\lambda^*_i(dt)Eh_i(\theta)\textbf{\emph{1}}_{\{\theta\leq
T\}},~~t\in\Lambda\cap[0,T],\ i=1,2.
\end{equation}
That is to say, $\tilde{\lambda}_i(t),\ i=1,2,$ are proportional to the renewal
functions of some corresponding delayed renewal counting processes.

Similarly, we define a r.v.
$\theta^{**}_1$ with a proper distribution given by
\begin{equation}\label{rem2110}
P(\theta^{**}_1\in dt)=\frac{g(t)}{Eg(\theta)\emph{\textbf{1}}_{\{\theta\leq
T\}}}G( dt),~~t\in\Lambda\cap[0,T],
\end{equation}
and a sequence of i.i.d. r.v.$^,$s $\{\theta^{**}_{k}:k\in\mathbb{N}\}$ with the same distribution as $\theta_1^{**}$. Then the
inter-arrival times $\theta^{**}_{1},\theta_k, k\ge 2,$ also
follow a delayed renewal counting process $\{N^{**}(t); t\geq0\}$
with a corresponding renewal function $\lambda^{**}(t)=EN^{**}(t)$,
$t\geq0$, such that for any $t\in\Lambda\cap[0,T]$
\begin{eqnarray}\label{rem21200}
\tilde{\tilde{\lambda}}(dt)=\lambda^{**}(dt)Eg(\theta)\textbf{\emph{1}}_{\{\theta\leq T\}}\ge P(\theta^{**}_1\in dt)Eg(\theta)\textbf{\emph{1}}_{\{\theta\leq
T\}}=g(t)G(dt).
\end{eqnarray}
\end{remark}


\begin{remark}
If $\lim_{z\rightarrow\infty}g_{i,j}(z,s)=g_i(s)$ holds uniformly for all $s>0$ and $1\le i\neq j\le2$, then Condition 2 is automatically satisfied given Conditions 1 and 3.
\end{remark}

The remaining of this paper consists of four sections. Section 2 gives several lemmas which are pivotal to the proof of our main result in Section 3. And Section 4 presents some concrete joint distributions or copulas for the dependent structure among the claim sizes and the inter-arrival time to demonstrate the results we established. Finally, a non locally almost decreased distribution in the class $\mathcal{S}_{loc}$ is given in Section 5.

\setcounter{equation}{0}
\section{Some lemmas}

This section collects three technical lemmas
having their own independent interests. The first lemma slightly improved a method in Asmussen et al. \cite{AFK2003}, where it was requested that $g(x)x^{-1}\to 0$.

\begin{lemma}\label{l32}
If $F\in\mathcal{L}_{\Delta_d}$ for some $0<d\le\infty$, then there exists a function $g(\cdot): [0,\infty)\mapsto(0,\infty)$ such
that $g(x)\uparrow\infty,\ g(x)x^{-1}\downarrow 0$ and
\begin{equation}\label{lm4200}
\overline{F}(x+y+\Delta_d)\sim\overline{F}(x+\Delta_d)\ uniformly\ for\ all\ |y|\leq g(x).
\end{equation}
\end{lemma}

{\it{Proof.}}\ \ Since $F\in\mathcal{L}_{\Delta_d}$ with some $0<d\le\infty$, for each integer $n\ge1$ and all $\mid y\mid\le n$, there exists a number $x_n>\max\{2x_{n-1},n^2\}$ such that, when $x\ge x_n$,
$$(1-n^{-1})F(x+\Delta_d)<F(x+y+\Delta_d)<(1+n^{-1})F(x+\Delta_d),$$
where $x_0=1$. Let $g_0(\cdot): [0,\infty)\mapsto(0,\infty)$ be a function such that
$$g_0(x)=\textbf{1}(0\le x < 1)+\sum_{n=1}^\infty n\textbf{1}(x_n\le x < x_{n+1}).$$
Clearly, $g_0(x)\uparrow\infty,\ g_0(x)x^{-1}\rightarrow 0$ and $\overline{F}(x+y+\Delta_d)\sim\overline{F}(x+\Delta_d)$
uniformly for all $|y|\leq g_0(x)$. Define a continuous linear function $g(\cdot): [0,\infty)\mapsto(0,\infty)$ by
$$g(x)=\textbf{1}(0\le x < 1)+\sum_{n=1}^\infty \Big((x_{n+1}-x_n)^{-1}x+\big(n-x_n(x_{n+1}-x_n)^{-1}\big)\Big)\textbf{1}(x_n\le x < x_{n+1}).$$
Then $g_0(x)\ge g(x)\uparrow\infty$ for all $x\ge0$, and when $x_n\le x\le x_{n+1}$,
$$g_0(x)x^{-1}=(x_{n+1}-x_n)^{-1}+\big(n-x_n(x_{n+1}-x_n)^{-1}\big)x^{-1}\downarrow0,$$
further (\ref{lm4200}) holds. \hfill$\Box$
\begin{lemma}\label{l31}
A distribution $F\in\mathcal{L}_{loc}$ if and only if for any $g>0$ and any $b>a>0$,
\begin{equation}\label{lm420}
\sup_{ |y|\le g,\ d,s\in(a,b]} |F(x+y+\Delta_d)\big(F(x+\Delta_s)\big)^{-1}-ds^{-1}|\to 0.
\end{equation}
\end{lemma}

{\it{Proof.}}\ \ We only need to prove (\ref{lm420}) with $s=1$ by $F\in\mathcal{L}_{loc}$. For any fixed integer $n\ge2$, we denote $d_k=a+k(b-a)n^{-1}$ for all integers $1\le k\le n$. From Lemma 2.1 of Wang and Wang (2011) and uniformly convergent theorem of slow variable function, we know that for all fixed $h>0$,
\begin{equation}\label{lm422}
\sup_{ |y|\le h,\ 1\le k\le n} |F(x+y+\Delta_{d_k})\big(F(x+y+\Delta_1)\big)^{-1}-d_k|\to 0.
\end{equation}
For any $d\in(a,b]$, there is an integer $1\le k\le n$ such that $d\in(d_{k-1},d_k]$. For any $0<\varepsilon<1$, we take $n$ large enough such that
$$\sqrt{1-\varepsilon}d<d_{k-1}<d_k<\sqrt{1+\varepsilon}d.$$
Then by (\ref{lm422}), there is a constant $x_0=x_0(F,\varepsilon, h, d_k, 1\le k\le n)>0$ such that, for all $x\ge x_0$ and $\mid y\mid\le h$,
\begin{eqnarray*}
F(x+y+\Delta_d)\big(F(x+y+\Delta_1)\big)^{-1}&\le& F(x+y+\Delta_{d_k})\big(F(x+y+\Delta_1)\big)^{-1}\nonumber\\
&\le&\sqrt{1+\varepsilon}d_k\le (1+\varepsilon)d
\end{eqnarray*}
and
\begin{eqnarray*}
F(x+y+\Delta_d)\big(F(x+y+\Delta_1)\big)^{-1}&\ge&F(x+y+\Delta_{d_k})\big(F(x+y+\Delta_1)\big)^{-1}\nonumber\\
&\ge&\sqrt{1-\varepsilon}d_{k-1}\ge (1-\varepsilon)d.
\end{eqnarray*}
Therefore, according to arbitrary of $\varepsilon$, (\ref{lm420}) holds.\hfill$\Box$

\begin{lemma}\label{l33}
Consider the bidimensional risk model (\ref{p2})
satisfying Conditions 1, 3 and 4.
If $F_i\in\mathcal{S}_{loc}$, $i=1,2$, then for any $T\in\Lambda$ and every fixed
$n\geq1$, it holds uniformly for all $t\in\Lambda\cap(0,T]$ that
\begin{eqnarray}\label{lm43}
&&P\Big(\sum\limits_{k=1}^{n}\vec{X}_ke^{-r\sigma_k}\in\vec{x}+\vec{\Delta},N(t)=n\Big)\nonumber\\
&\sim&\sum\limits_{k=1}^{n}\sum_{j=1}^{n} P(X^{(1)}_ke^{-r\sigma_k}\in x_1+\Delta^{(1)},X^{(2)}_je^{-r\sigma_j}\in x_2+\Delta^{(2)},N(t)=n).
\end{eqnarray}
\end{lemma}
\vspace{0.2cm}

{\it{Proof.}}\ \ For every $t\in\Lambda\cap[0,T]$ and $n\geq1$, we write
$$\Omega_{n}(t)=\Big\{(s_1,\ldots,s_n)\in[0,t]^n: t_n=\sum_{k=1}^ns_k\leq t\Big\}.$$

Firstly, for every fixed $n\geq1$, we prove that it holds uniformly for all $(s_1,\cdots,s_n)\in\Omega_n(t),z_l\ge0,1\le l\le n$, and $t\in\Lambda\cap[0,T]$ that
\begin{eqnarray}\label{303}
&&P\Big(\sum_{k=1}^{n}X^{(1)}_ke^{-rt_k}\in x_1+\Delta^{(1)}|X^{(2)}_l=z_l,\theta_l=s_l,1\le l\le n\Big)\nonumber\\
&\sim&\sum_{k=1}^{n}P(X^{(1)}_ke^{-rt_k}\in x_1+\Delta^{(1)}|X^{(2)}_k=z_k,\theta_k=s_k).
\end{eqnarray}
Let us proceed by induction. Clearly, the assertion holds for
$n=1$. Now we assume that the assertion holds for $n=m-1$.
Then when $n=m$, due to the fact $F_1\in\mathcal{S}_{loc}\subset\mathcal{L}_{loc}$, by Lemma \ref{l32}, there exists a function $g_1(\cdot): [0,\infty)\mapsto(0,\infty)$ such
that $g_1(x)\rightarrow\infty$, $g_1(x)x^{-1}\downarrow0$ and
$\overline{F_1}(x_1+y+\Delta^{(1)})\sim\overline{F_1}(x_1+\Delta^{(1)})$ holds uniformly for all
$|y|\leq (g_1(x)+2d_1)e^{rT}$. Thus, by standard methods, we have
\begin{eqnarray}
&&P\Big(\sum_{k=1}^{m}X^{(1)}_ke^{-rt_k}\in x_1+\Delta^{(1)}|X^{(2)}_l=z_l,\theta_l=s_l,1\le l\le m\Big)\nonumber\\
&=&P\Big(\sum_{k=1}^{m}X^{(1)}_ke^{-rt_k}\in x_1+\Delta^{(1)},\sum_{k=1}^{m-1}
X^{(1)}_ke^{-rt_k}\leq g_1(x_1)|X^{(2)}_l=z_l,\theta_l=s_l,1\le l\le m\Big)\nonumber\\
&&\ \ \ \ \ \ \ \ \ \ +P\Big(\sum_{k=1}^{m}X^{(1)}_ke^{-rt_k}\in x_1+\Delta^{(1)},\sum_{k=1}^{m-1}X^{(1)}_ke^{-rt_k}>x_1-g_1(x_1)\nonumber\\
&&\ \ \ \ \ \ \ \ \ \ \ \ \ \ \ \ \ \ \ \ |X^{(2)}_l=z_l,\theta_l=s_l,1\le l\le m\Big)\nonumber\\
&&\ \ \ \ \ \ \ \ \ \ +P\Big(\sum_{k=1}^{m}X^{(1)}_ke^{-rt_k}\in x_1+\Delta^{(1)},g_1(x_1)
<\sum_{k=1}^{m-1}X^{(1)}_ke^{-rt_k}\leq x_1-g_1(x_1)\nonumber\\
&&\ \ \ \ \ \ \ \ \ \ \ \ \ \ \ \ \ \ \ \ |X^{(2)}_l=z_l,\theta_l=s_l,1\le l\le m\Big)\nonumber\\
&=&I_1(x_1,m)+I_{2}(x_1,m)+I_{3}(x_1,m).\label{304}
\end{eqnarray}

For $I_{1}(x_1,m)$, by $F_1\in\mathcal {S}_{loc}\subset\mathcal {L}_{loc}$, (\ref{207})
and Lemma \ref{l31}, it holds uniformly for all $(s_1,\ldots,s_m)\in\Omega_m(t),z_l\ge0,1\le l\le n$,
and $t\in\Lambda\cap[0,T]$ that
\begin{eqnarray}
I_{1}(x_1,m)&=&\int_0^{g_1(x_1)} P(X^{(1)}_me^{-rt_m}\in x_1-y+\Delta^{(1)}|X^{(2)}_m=z_m,\theta_m=s_m)\nonumber\\
&&\ \ \ \ \ \ \ \ \ \ \ \cdot P(\sum_{k=1}^{m-1}X^{(1)}_ke^{-rt_k}\in dy|X^{(2)}_l=z_l,\theta_l=s_l,1\le l\le m-1)\nonumber\\
&\sim& P(X^{(1)}_me^{-rt_m}\in x_1+\Delta^{(1)}|X^{(2)}_m=z_m,\theta_m=s_m).\label{305}
\end{eqnarray}

For $I_2(x_2,m)$, by the induction assumption, $F_1\in\mathcal {S}_{loc}$, (\ref{207}) and Lemma \ref{l31},  it holds uniformly for all
$(s_1,\ldots,s_m)\in\Omega_m(t),z_l\ge0,1\le l\le n$ and $t\in\Lambda\cap[0,T]$ that
\begin{eqnarray}
&&I_{2}(x_1,m)=P\Big(\sum\limits_{k=1}^{m}X^{(1)}_ke^{-rt_k}\in x_1+\Delta^{(1)},X^{(1)}_me^{-rt_m}\le g_1(x_1)+d_1,\nonumber\\
&&\ \ \ \ \ \ \ \ \ \ \ \ \ \ \ \ \ \ \ \ x_1-g_1(x_1)<\sum\limits_{k=1}^{m-1}X^{(1)}_ke^{-rt_k}
\le x_1+d_1|X^{(2)}_l=z_l,\theta_l=s_l,1\le l\le m\Big)\nonumber\\
&\leq&\int_0^{g_1(x_1)+d_1}P\Big(\sum\limits_{k=1}^{m-1}X^{(1)}_ke^{-rt_k}
\in x_1-y+\Delta^{(1)}|X^{(2)}_l=z_l,\theta_l=s_l,1\le l\le m-1\Big)\nonumber\\
&&\ \ \ \ \ \ \ \ \ \ \ \ \ \ \ \ \ \ \ \ \cdot P(X^{(1)}_me^{-rt_m}\in dy\mid X^{(2)}_m=z_m,\theta_m=s_m)\nonumber\\
&\sim&\sum\limits_{k=1}^{m-1}P\Big(X^{(1)}_ke^{-rt_k}
\in x_1+\Delta^{(1)}|X^{(2)}_k=z_k,\theta_k=s_k\Big)\label{306}
\end{eqnarray}
and
\begin{eqnarray}
&&I_{2}(x_1,m)\ge P\Big(\sum\limits_{k=1}^{m}X^{(1)}_ke^{-rt_k}\in x_1+\Delta^{(1)},X^{(1)}_me^{-rt_m}\le g_1(x_1),\nonumber\\
&&\ \ \ \ \ \ \ \ \ \ \ \ \ \ \ \ \ \ \ \ x_1-g_1(x_1)<\sum\limits_{k=1}^{m-1}X^{(1)}_ke^{-rt_k}
\le x_1+d_1|X^{(2)}_l=z_l,\theta_l=s_l,1\le l\le m\Big)\nonumber\\
&=&P\Big(\sum\limits_{k=1}^{m}X^{(1)}_ke^{-rt_k}\in x_1+\Delta^{(1)},X^{(1)}_me^{-rt_m}\le g_1(x_1)
|X^{(2)}_l=z_l,\theta_l=s_l,1\le l\le m\Big)\nonumber\\
&\sim& \sum\limits_{k=1}^{m-1}P\Big(X^{(2)}_ke^{-rt_k}
\in x_1+\Delta^{(1)}|X^{(2)}_k=z_k,\theta_k=s_k\Big).\label{3066}
\end{eqnarray}

Now, we analyze $I_3(x_1,m)$. We define two positive functions $g_1^*(\cdot)$ and $g_1^{**}(\cdot)$ by $g_1^*(\cdot)=\max\{0,g_1(\cdot)-d_1\}$ and $g_1^{**}(\cdot)=\max\{0,g_1^*(\cdot)-d_1\}$, respectively. By the induction assumption, $F_1\in\mathcal {S}_{loc}$, (\ref{207}), Lemma \ref{l31} and Condition 4, it holds uniformly for all $(s_1,\ldots,s_m)\in\Omega_m(t),z_l\ge0,1\le l\le n,$ and
$t\in\Lambda\cap[0,T]$ that
\begin{eqnarray*}
&&I_{3}(x_1,m)=P\Big(\sum_{k=1}^{m}X^{(1)}_ke^{-rt_k}\in x_1+\Delta^{(1)},g_1(x_1)
<\sum_{k=1}^{m-1}X^{(1)}_ke^{-rt_k}\leq x_1-g_1(x_1),\nonumber\\
&&\ \ \ \ \ \ \ \ \ \ \ \ \ \ \ \ \ \ \ \ g_1(x_1)<X^{(1)}_me^{-rt_m}\le x_1-g_1(x_1)+d_1|X^{(2)}_l=z_l,\theta_l=s_l,1\le l\le m\Big)\nonumber\\
&\le&\int_{g_1^*(x_1)}^{x_1-g_1^*(x_1)}P\Big(\sum_{k=1}^{m-1}X^{(1)}_ke^{-rt_k}\in x_1-y+\Delta^{(1)}|X^{(2)}_l=z_l,\theta_l=s_l,1\le l\le m-1\Big)\nonumber\\
&&\ \ \ \ \ \ \ \ \ \ \ \ \ \ \ \ \ \ \ \cdot P(X^{(1)}_me^{-rt_m}\in dy|X^{(2)}_m=z_m,\theta_m=s_m)\nonumber\\
&\sim&\sum_{k=1}^{m-1}\int_{g_1^*(x_1)}^{x_1-g_1^*(x_1)}P(X^{(1)}_ke^{-rt_k}\in x_1-y+\Delta^{(1)}|X^{(2)}_k=z_k,\theta_k=s_l)\nonumber\\
&&\ \ \ \ \ \ \ \ \ \ \ \ \ \ \ \ \ \ \ \cdot P(X^{(1)}_me^{-rt_m}\in dy|X^{(2)}_m=z_m,\theta_m=s_m)\nonumber\\
&=&O\Big(\sum_{k=1}^{m-1}\int_{g_1^*(x_1)}^{x_1-g_1^*(x_1)}P(X^{(1)}_ke^{-rt_k}\in x_1-y+\Delta^{(1)})\nonumber\\
&&\ \ \ \ \ \ \ \ \ \ \ \ \ \ \ \ \ \ \ \cdot P(X^{(1)}_me^{-rt_m}\in dy|X^{(2)}_m=z_m,\theta_m=s_m)\Big)\nonumber\\
&=&O\Big(\sum_{k=1}^{m-1}P\big(X^{(1)}_ke^{-rt_k}+X^{(1)}_me^{-rt_m}\in x_1+\Delta^{(1)},X^{(1)}_ke^{-rt_k}\in(g_1^*(x_1),x_1-g_1^*(x_1)+d_1],\nonumber\\
&&\ \ \ \ \ \ \ \ \ \ \ \ \ \ \ \ \ \ \ X^{(1)}_me^{-rt_m}\in(g_1^*(x_1),x_1-g_1^*(x_1)]|X^{(2)}_m=z_m,\theta_m=s_m\big)\Big)\nonumber\\
&=&O\Big(\sum_{k=1}^{m-1}\int_{g_1^{**}(x_1)}^{x_1-g_1^{**}(x_1)}P(X^{(1)}_me^{-rt_m}\in x_1-y+\Delta^{(1)})P(X^{(1)}_ke^{-rt_k}\in dy)\Big)\nonumber\\
&=&O\Big(\sum_{k=1}^{m-1}\int_{g_1^{**}(x_1)}^{x_1-g_1^{**}(x_1)}P(X^{(1)}_ke^{-rt_k}\in x_1-y+\Delta^{(1)})P(X^{(1)}_ke^{-rt_k}\in dy)\Big)\nonumber\\
&=&O\Big(\sum_{k=1}^{m-1}\int_{e^{rt_k}g_1^{**}(x_1)}^{e^{rt_k}(x_1-g_1^{**}(x_1))}P(X^{(1)}_k\in e^{rt_k}x_1-u+\Delta^{(1)})P(X^{(1)}_k\in du)\Big).
\end{eqnarray*}
Since $x^{-1}g_1(x)\downarrow0$, $ag_1(x)\ge g_1(ax)$ for any $a>1$. Thus, by $F_1\in\mathcal{S}_{loc}$, (\ref{206}) and Lemma \ref{l31}, it holds uniformly for all $(s_1,\ldots,s_m)\in\Omega_m(t),z_l\ge0,1\le l\le n,$ and
$t\in\Lambda\cap[0,T]$ that
\begin{eqnarray}
I_{3}(x_1,m)&=&O\Big(\sum_{k=1}^{m-1}\int_{g_1^{**}(e^{rt_k}x_1)}^{e^{rt_k}x_1-g_1^{**}(e^{rt_k}x_1)}P(X^{(1)}_k\in e^{rt_k}x_1-u+\Delta^{(1)})P(X^{(1)}_k\in du)\Big)\nonumber\\
&=&o\Big(\sum_{k=1}^{m-1}P(X^{(1)}_ke^{-rt_k}\in x_1+\Delta^{(1)})\Big)\nonumber\\
&=&o\Big(\sum_{k=1}^{m}P(X^{(1)}_ke^{-rt_k}\in x_1+\Delta^{(1)}|X^{(2)}_k=y_k,\theta_k=s_k)\Big).\label{307}
\end{eqnarray}

Therefore, combining (\ref{304})-(\ref{307}) implies that (\ref{303}) holds uniformly for all $(s_1,\ldots,s_m)\in\Omega_m(t)$, $z_l\ge0,1\le l\le n,$ and $t\in\Lambda\cap[0,T]$.

Secondly, following the same manner, we can similarly apply Conditions 3, 4 and Lemma \ref{l31} to obtain the following result: For every fixed $1\le k\le n$, it holds for all $(s_1,\ldots,s_n)\in\Omega_n(t),y_k\ge0,$
and $t\in\Lambda\cap[0,T]$ that
\begin{eqnarray}\label{3030}
&&P\Big(\sum_{j=1}^{n}X^{(2)}_je^{-rt_j}\in x_2+\Delta^{(2)}|X^{(1)}_k=y_k,\theta_l=s_l,1\le l\le n\Big)\nonumber\\
&\sim&\sum_{j=1}^{n}P(X^{(2)}_je^{-rt_j}\in x_2+\Delta^{(2)}|X^{(1)}_k=y_k,\theta_j=s_j).
\end{eqnarray}

Finally, for every fixed $n\geq1$, by (\ref{303}) and (\ref{3030}), it holds uniformly
for all $t\in\Lambda\cap[0,T]$ that
\begin{eqnarray*}
&&P\Big(\sum\limits_{k=1}^{n}\vec{X}_ke^{-r\sigma_k}\in\vec{x}+\vec{\Delta},
N(t)=n\Big)\nonumber\\
&\sim&\idotsint\limits_{\Omega_{n}(t)}\idotsint\limits_{\{\sum_{j=1}^{n}z_je^{-rt_j}\in x_2+\Delta^{(2)}\}}
\sum\limits_{k=1}^{n}P\Big(X^{(1)}_ke^{-rt_k}\in x_1+\Delta^{(1)}|X_k^{(2)}=z_k,\theta_k=s_k\Big)\nonumber\\
&&\ \ \ \ \ \ \ \ \ \ \ \ \ \ \ \ \ \ \ \ \ \ \ \ \cdot\prod_{k=1}^nP(X_k^{(2)}\in dz_k|\theta_k=s_k)\overline{G}(t-t_n)\prod_{l=1}^nG(ds_l)\nonumber\\
&\sim&\sum_{k=1}^{n}\sum_{j=1}^{n}\idotsint\limits_{\Omega_{n}(t)}\int_{x_1
e^{rt_k}}^{(x_1+d_1)e^{rt_k}}
P\Big(X^{(2)}_j e^{-rt_j}\in x_2+\Delta^{(2)}|X^{(1)}_k=y_k,\theta_j=s_j\Big)\nonumber\\
&&\ \ \ \ \ \ \ \ \ \ \ \ \ \ \ \ \ \ \ \ \ \ \ \ \cdot P(X^{(1)}_k\in dy_k\mid\theta_k=s_k)\overline{G}(t-t_n)\prod_{l=1}^nG(ds_l)\nonumber\\
&=&\sum\limits_{k=1}^{n}\sum_{j=1}^{n} P\Big(X^{(1)}_ke^{-r\sigma_k}\in x_1+\Delta^{(1)},X^{(2)}_je^{-r\sigma_j}\in x_2+\Delta^{(2)},N(t)=n\Big).
\end{eqnarray*}
Thus, the proof is completed.  \hfill$\Box$

\vspace{0.2cm}

\setcounter{equation}{0}
\setcounter{lemma}{0}
\section{Proof of Theorem 2.1}
\quad~~

For any fixed integer $N\geq 1$,
\begin{align}
F_{\vec{D}_{r}(t)}(\vec{x}+\vec{\Delta})&=\Big(\sum_{n=1}^N+\sum_{n=N+1}^{\infty}\Big)
P\Big(\sum_{k=1}^n\vec{X}_{k}e^{-r\sigma_k}\in\vec{x}+\vec{\Delta},N(t)=n\Big)\nonumber\\
&=J_1(\vec{x},t)+J_2(\vec{x},t).\label{5001}
\end{align}
By Lemma \ref{l33}, it holds
uniformly for all $t\in \Lambda\cap[0,T]$ that
\begin{eqnarray}
J_1(\vec{x},t)&\sim&\Bigg(\sum_{n=1}^{\infty}-\sum_{n=N+1}^{\infty}\Bigg)\sum\limits_{k=1}^{n}\sum_{j=1}^{n} P\Big(X^{(1)}_ke^{-r\sigma_k}\in x_1+\Delta^{(1)},\nonumber\\
&&\ \ \ \ \ \ \ \ \ \ \ \ X^{(2)}_je^{-r\sigma_j}\in x_2+\Delta^{(2)}, N(t)=n\Big)\nonumber\\
&=&J_{11}(\vec{x},t)-J_{12}(\vec{x},t).\label{5002}
\end{eqnarray}
Moreover,
\begin{eqnarray}
J_{11}(\vec{x},t)&=&\sum_{k=1}^{\infty}\sum_{j=1}^{k-1} P(X^{(1)}_ke^{-r\sigma_k}\in x_1+\Delta^{(1)},X^{(2)}_je^{-r\sigma_j}\in x_2+\Delta^{(2)},\sigma_k\leq t)\nonumber\\
&&\ \ \ +\sum_{j=1}^{\infty}\sum_{k=1}^{j-1} P(X^{(1)}_ke^{-r\sigma_k}\in x_1+\Delta^{(1)},X^{(2)}_je^{-r\sigma_j}\in x_2+\Delta^{(2)},\sigma_j\leq t)\nonumber\\
&&\ \ \ \ \ \ \ \ \ \ +\sum_{k=1}^{\infty}P(\vec{X}_{k}e^{-r\sigma_k}\in\vec{x}+\vec{\Delta}, \sigma_k\leq t)\nonumber\\
&=&\sum_{i=1}^3 J_{11i}(\vec{x},t),\label{5031}
\end{eqnarray}
where $\sum_{k=1}^0=\sum_{j=1}^0=0$.

Denote $c_1=Eh_1(\theta)\textbf{1}_{\{\theta\leq
T\}}Eh_2(\theta)\textbf{1}_{\{\theta\leq T\}}$. According to (\ref{202}), (\ref{rem211}) and
(\ref{rem212}), and noting that
$\sigma_{k-1}-\sigma_j+\theta_{1,k}^*,\sigma_{k-j-1}+\theta_{1,k-j}^*$ and
$\theta_{1,1}^*+\sigma_{k-j}-\theta_{1}$ are identically distributed
r.v.s for $k\ge j+1$ and $j\ge 1$, it holds uniformly for all $t\in
\Lambda\cap[0,T]$ that
\begin{eqnarray}
&&J_{111}(\vec{x},t)=\sum_{k=1}^{\infty}\sum_{j=1}^{k-1}\int_{0-}^t\int_{0-}^{t-u}\int_{0-}^{t-u-v}
\int_{0-}^{t-u-v-w}P(X^{(1)}_{k}e^{-r(u+v+w+z)}\in x_1+\Delta^{(1)} \nonumber\\
&&\ \ \ \ \ \ \ \ \ \ \ \ \ \ \ \ \ \ |\theta_k=u)P(X^{(2)}_{j}e^{-r(w+z)}\in x_2+\Delta^{(2)}|\theta_j=w)\nonumber\\
&&\ \ \ \ \ \ \ \ \ \ \ \ \ \ \ \ \ \ \ \ \ \ \ \ \ \ \ \ \ \ \cdot P(\sigma_{j-1}\in dz)P(\theta_j\in dw)P(\sigma_{k-1}-\sigma_j\in dv)P(\theta_k\in du)\nonumber\\
&\sim&\sum_{k=1}^{\infty}\sum_{j=1}^{k-1}\int_{0-}^t\int_{0-}^{t-u}\int_{0-}^{t-u-v}\int_{0-}^{t-u-v-w}
P(X^{(1)}_{k}e^{-r(u+v+w+z)}\in x_1+\Delta^{(1)})h_1(u)\nonumber\\
&&\ \ \ \ \ \ \ \ \ \ \ \ \ \ \ \ \ \ \ \ \cdot P( X^{(2)}_{j}e^{-r(w+z)}\in x_2+\Delta^{(2)})h_2(w)\nonumber\\
&&\ \ \ \ \ \ \ \ \ \ \ \ \ \ \ \ \ \ \ \ \ \ \ \ \ \ \ \ \ \ \cdot P(\sigma_{j-1}\in dz)P(\theta_j\in dw)P(\sigma_{k-1}-\sigma_j\in dv)P(\theta_k\in du)\nonumber\\
&=&c_1\int_{0-}^t\int_{0-}^{t-u}P(X^{(1)}e^{-r(u+v)}\in(x_1,x_1+d_1])P( X^{(2)}e^{-rv}\in(x_2,x_2+d_2])
\nonumber\\
&&\ \ \ \ \ \ \ \ \ \ \ \ \ \ \ \ \ \ \ \ \cdot\sum_{j=1}^{\infty}\Big(\sum_{k=j+1}^{\infty}P(\sigma_{k-1}-\sigma_j+\theta_{1,k}^*\in du)\Big)P(\sigma_{j-1}+\theta_{2,j}^*\in dv)\nonumber\\
&=&c_1\int_{0-}^t\int_{0-}^{t-u}P(X^{(1)}e^{-r(u+v)}\in x_1+\Delta^{(1)})P( X^{(2)}e^{-rv}\in x_2+\Delta^{(2)})\nonumber\\
&&\ \ \ \ \ \ \ \ \ \ \ \ \ \ \ \ \ \ \ \ \cdot\sum_{j=1}^{\infty}\Big(\sum_{l=1}^{\infty}P(\sigma_{l-1}+\theta_{1,l}^*\in du)\Big)
P(\sigma_{j-1}+\theta_{2,j}^*\in dv)\nonumber\\
&=&\int_{0-}^{t}\int_{0-}^{t-u}P(X^{(1)}e^{-r(u+v)}\in x_1+\Delta^{(1)})P( X^{(2)}e^{-rv}\in x_2+\Delta^{(2)})
\tilde{\lambda}_2(dv)\tilde{\lambda}_1(du). \label{5034}
\end{eqnarray}
Similarly, it holds uniformly for all $t\in\Lambda\cap[0,T]$ that
\begin{eqnarray}
&&J_{112}(\vec{x},t)\sim\int_{0-}^{t}\int_{0-}^{t-u}P(X^{(1)}e^{-ru}\in x_1+\Delta^{(1)})P( X^{(2)}e^{-r(u+v)}\in x_2+\Delta^{(2)})\nonumber\\
&&\ \ \ \ \ \ \ \ \ \ \ \ \ \ \ \ \ \ \ \ \cdot\tilde{\lambda}_2(dv)\tilde{\lambda}_1(du).\label{5035}
\end{eqnarray}
Now, we shall analyze $J_{113}(\vec{x},t)$. Denote
$c_2=Eg(\theta)\textbf{1}_{\{\theta\leq T\}}$. According to Condition 2
and (\ref{rem2110}), it holds uniformly for all $t\in\Lambda\cap[0,T]$
\begin{eqnarray}
&&J_{113}(\vec{x},t)\sim\sum_{k=1}^{\infty}\int_{0-}^t\int_{0-}^{t-u}\prod_{i=1}^2 P\big(X^{(i)}_ke^{-ru}\in x_i+\Delta^{(i)})
g(u)P(\sigma_{k-1}\in dv\big)G(du)\nonumber\\
&=&c_2\sum_{k=1}^{\infty}\int_{0-}^t
\prod_{i=1}^2 P\big(X^{(i)}_ke^{-ru}\in x_i+\Delta^{(i)})P(\sigma_{k-1}+\theta^{**}_{1,k}\in du\big)\nonumber\\
&=&\int_{0-}^{t}\prod_{i=1}^2 P(X^{(i)}_ke^{-ru}\in x_i+\Delta^{(i)})\Big)
\tilde{\tilde{\lambda}}(du).\label{503}
\end{eqnarray}

Next, we focus on the analysis of $J_{12}(\vec{x},t)$ and denote it by
\begin{eqnarray}
&&J_{12}(\vec{x},t)=\sum_{n=N+1}^{\infty}\sum_{1\le k,j\le n}\int_{0-}^{t}\cdot\cdot\cdot\int_{0-}^{t-\sum_{i=1}^{n-1}u_i}P(X^{(1)}_ke^{-r\sum_{i=1}^{k}u_i}\in x_1+\Delta^{(1)}, \nonumber\\
&&\ \ \ \ \ X^{(2)}_j e^{-r\sum_{i=1}^{j}u_i}\in x_2+\Delta^{(2)}|\theta_k=u_k,\theta_j=u_j)P(N(t-u_1)>n-1)\prod_{i=1}^nG(du_i)\nonumber\\
&=&\sum_{n=N+1}^{\infty}\sum_{1\le k,j\le n}\int_{0-}^{t}\cdot\cdot\cdot\int_{0-}^{t-\sum_{i=1}^{n-1}u_i}P_{k,j}(\vec{x},\vec{u_1})P(N(t-u_1)>n-1)\prod_{i=1}^nG(du_i).\label{5041}
\end{eqnarray}
We first deal with $P_{k,j}(\vec{x},\vec{u})$
for two cases that $k=j$ and $k\neq j$, respectively.
In every case above, by (\ref{201}) or (\ref{206}), then by Condition 4, it holds uniformly for all $u_i\in[0,t],1\le i\le n,$ and $t\in \Lambda\cap[0,T]$ that
\begin{align}
P_{k,j}(\vec{x},\vec{u})&=O\Big(P(X^{(1)}_ke^{-ru_1}\in x_1+\Delta^{(1)})P(X^{(2)}_j e^{-ru_1}\in x_2+\Delta^{(2)})\Big).\label{50410}
\end{align}
From (\ref{5041}), (\ref{50410}), (\ref{rem2110}) and (\ref{rem21200}), it follows that
\begin{align}
&J_{12}(\vec{x},t)=O\Big(\int_{0-}^t\prod_{i=1}^2 P(X^{(i)}_ke^{-ru}\in x_i+\Delta^{(i)})G(du)
\sum_{n=N+1}^{\infty}n^2 P(N(T)>n-1)\Big)\nonumber\\
&=O\Big(\int_{0-}^t \prod_{i=1}^2 P(X^{(i)}_ke^{-ru}\in x_i+\Delta^{(i)})\tilde{\tilde{\lambda}}(du)
\sum_{n=N+1}^{\infty}n^2P(N(T)>n-1)\Big)\label{5042}
\end{align}
holds uniformly for all $t\in \Lambda\cap[0,T]$.
Then, combining (\ref{503}), (\ref{5042})
and $E(N(T))^3<\infty$ gives
\begin{eqnarray}
&&\lim_{N\rightarrow\infty}\limsup
\sup_{t\in\Lambda\cap[0,T]}\frac{J_{12}(\vec{x},t)}{J_{113}(\vec{x},t)}=0.\label{505}
\end{eqnarray}

Finally, we deal with $J_{2}(\vec{x},t)$.
By (\ref{2011}), Condition 4, Kesten's inequality and Lemma \ref{l31} for some $\epsilon>0$ satisfying $(1+\alpha)(1+\epsilon)<1+\beta$, we have
\begin{eqnarray*}
&&J_2(x,t)\le (1+C_2)\sum_{n=N+1}^\infty\int\cdot\cdot\cdot\int_{0\le s_k\le t-\sum_{i=1}^{k-1}s_i,1\le k\le n}\nonumber\\
&&\ \ \ \ \ \ \ \ \ \ \ \int\cdot\cdot\cdot\int_{\sum_{k=1}^nu_k e^{-r\sum_{i=1}^{k} s_i}\in x_1+\Delta^{(1)}}\int\cdot\cdot\cdot\int_{\sum_{j=2}^nv_je^{-r\sum_{i=1}^{j} s_i}\in x_2+\Delta^{(2)}}\nonumber\\
&&\ \ \ \ \ \ \ \ \ \ \ \ \ \ \ \ \ \ \ \prod_{k=1}^nP(X_k^{(1)}\in du_k)P(X_k^{(2)}\in dv_k)P(\theta_k\in ds_k)\nonumber\\
&\le& (1+C_2)\sum_{n=N+1}^\infty\int\cdot\cdot\cdot\int_{0\le s_k\le t-\sum_{i=1}^{k-1}s_i,1\le k\le n}\int\cdot\cdot\cdot\int_{\sum_{k=2}^nu_k e^{-r\sum_{i=1}^{k} s_i}\in x_1+\Delta^{(1)}}\nonumber\\
&&\ \ \ \ \ \ \ \ \ \ \ P(X_1^{(1)}e^{-rs_1}\in x_1-\sum_{k=2}^nu_k e^{r\sum_{i=1}^{k} s_i}+\Delta^{(1)})\int\cdot\cdot\cdot\int_{\sum_{j=2}^nv_je^{-r\sum_{i=2}^{j} s_i}\in x_2+\Delta^{(2)}}\nonumber\\
&&\ \ \ \ \ \ \ \ \ \ \ \ \ \ \ \ \ \ \ P(X_1^{(2)}e^{-rs_1}\in x_2-\sum_{k=2}^nv_k e^{-r\sum_{i=1}^{k} s_i}+\Delta^{(2)})\prod_{k=1}^nF_1(du_k)F_2(dv_k)G(ds_k)\nonumber\\
&=&O\Big(\sum_{n=N+1}^\infty\int\cdot\cdot\cdot\int_{0\le s_k\le t-\sum_{i=1}^{k-1}s_i,1\le k\le n}\int\cdot\cdot\cdot\int_{\sum_{k=2}^nu_k e^{r\sum_{i=1}^{k} s_i}\in x_1+\Delta^{(1)}}\nonumber\\
&&\ \ \ \ \ \ \ \ \ \ \ P(X_1^{(1)}e^{-rs_1}\in x_1-\sum_{k=2}^nu_k e^{-r s_1}+\Delta^{(1)})\int\cdot\cdot\cdot\int_{\sum_{j=2}^nv_je^{r\sum_{i=2}^{j} s_i}\in x_2+\Delta^{(2)}}\nonumber\\
&&\ \ \ \ \ \ \ \ \ \ \ \ \ \ \ \ \ \ \ P(X_1^{(2)}e^{-rs_1}\in x_2-\sum_{k=2}^nv_k e^{-r s_1}+\Delta^{(2)})\prod_{k=2}^nF_1(du_k)F_2(dv_k)G(ds_k)\Big)\nonumber\\
&=&O\Big(\sum_{n=N+1}^\infty\int_{0_-}^t F_1^{*n}((x_1+\Delta^{(1)})e^{rs})F_2^{*n}((x_2+\Delta_{(2)})e^{rs})G(ds)P(N(T)\ge n-1)\Big)\nonumber\\
&=&O\Big(\int_{0_-}^t F_1((x_1+\Delta^{(1)})e^{rs})F_2((x_2+\Delta_{(2)})e^{rs})G(ds)\sum_{n=N+1}^\infty(1+\alpha)^{n}P(N(T)\ge n-1)\Big)\nonumber\\
&=&O\Big(\int_{0_-}^t \prod_{i=1}^2 P(X^{(i)}_ke^{-ru}\in x_i+\Delta^{(i)})\tilde{\tilde{\lambda}}(du) Ee^{\beta N(T)}\textbf{1}_{\{ N(T)\ge N-1)\}}\Big)
\end{eqnarray*}
Thus,
\begin{equation}\label{5101}
\lim_{N\rightarrow\infty}\limsup
\sup_{t\in\Lambda\cap[0,T]}\frac{J_2(\vec{x},t)}{J_{113}(\vec{x},t)}=0.
\end{equation}
Combining (\ref{5001}), (\ref{5034}), (\ref{5035}), (\ref{503}), (\ref{505}) and (\ref{5101}), we obtain that (\ref{th211}) holds uniformly for all $t\in\Lambda\cap[0,T]$.

Therefore, (\ref{th212}) holds uniformly for all $t\in \Lambda\cap[0,T]$, following from
\begin{eqnarray*}\label{4012}
F_{-\vec{U}(\vec{x},t)}(\vec{\Delta})=P\big(\vec{U}(\vec{x},t)\in-\vec{\Delta}\big)
=F_{\vec{D}_r(t)}\Big(\vec{x}+\int_{0-}^te^{-ry}\vec{C}(dy)+e^{-rt}\vec{\Delta}\Big).
\end{eqnarray*}

\setcounter{equation}{0}
\setcounter{lemma}{0}
\section{Some copulas satisfying Conditions 1-3} \label{sec:copula}
\quad~~

The following concrete joint distributions or copulas of $(X^{(1)},X^{(2)},\theta)$ satisfying Conditions 1-3 show that the time-dependent structure in the Theorem \ref{thm201} has a larger range.
\vspace{0.2cm}

\indent{\bf Copula 4.1.}  Suppose that for all $k\ge 1,\ (X^{(1)}_k,X^{(2)}_k,\theta_k)$ follows a common tri-dimensional Sarmanov joint distribution
\begin{eqnarray*}
&&P(X^{(1)}_k\in dy,X^{(2)}_k\in dz,\theta_k\in ds)\nonumber\\
&=&\big(1+\gamma_{12}\phi_1(y)\phi_2(z)+\gamma_{13}\phi_1(y)\phi_3(s)+\gamma_{23}\phi_2(z)\phi_3(s)\big)
F_1(dy)F_2(dz)G(ds),
\end{eqnarray*}
where $\gamma_{ij}, 1\le i<j\le 3$, are constants and $\phi_i(\cdot),1\le i\le 3$, are continuous functions satisfying
$$1+\gamma_{12}\phi_1(y)\phi_2(z)+\gamma_{13}\phi_1(y)\phi_3(s)+\gamma_{23}\phi_2(z)\phi_3(s)\ge 0,\ y,z,s\in(-\infty,\infty)$$
and
$$E\phi_1(X^{(1)}_k)=E\phi_2(X^{(2)}_k)=E\phi_3(\theta_k)=0,$$
see Kotz et al. (2000). Therefore,
there exist constants $y_0,z_0$ and $s_0$ such that
\begin{eqnarray*}\label{2110}
\phi_1(y_0)=\phi_2(z_0)=\phi_3(s_0)=0.
\end{eqnarray*}

In particular, we have $\phi_1(y)=1-2F_1(y),\ y\in D_{X^{(1)}_1},\
\phi_2(z)=1-2F_2(z),\ z\in D_{X^{(2)}_1}$, $\phi_3(s)=1-2G(s),\
s\in D_{\theta_1}$ and $\gamma_{12},\gamma_{13},\gamma_{23}\in[-1,1]$, which give the well-known
tri-dimensional Farlie-Gumbel-Morgenstern (FGM) joint distribution,
where
$$D_X=\{x\in(-\infty,\infty):
P\big(X\in(x-\delta,x+\delta)\big)>0\ {\rm for\ all}\ \delta>0\}$$
for some random variable $X$.

In addition, we can also take $\phi_1(y)=e^{-y}-Ee^{-X^{(1)}_1}$,
$\phi_2(z)=e^{-z}-E e^{-X^{(2)}_1}$ and $\phi_3(s)=e^{-s}-E
e^{-\theta_1}$, or $\phi_1(y)=y^p-E({X^{(1)}_1})^p$,
$\phi_2(z)=z^p-E({X^{(2)}_1})^p$ and $\phi_3(s)=s^p-E\theta_1^p$,
and so on, for all $y\in D_{X^{(1)}_1}$, $z\in D_{X^{(2)}_1}$, $s\in
D_{\theta_1}$, where $p$ is some positive constant.

Furthermore, 
the distribution is required to satisfy the conditions that, for any $T\in\Lambda$,
\begin{eqnarray}\label{2111}
\inf_{y,z\in(-\infty,\infty),s\in[0,T]}\big(1+\gamma_{12}\phi_1(y)\phi_2(z)
+\gamma_{13}\phi_1(y)\phi_3(s)+\gamma_{23}\phi_2(z)\phi_3(s)\big)> 0
\end{eqnarray}
and for $i=1,2$, there exist positive constants $\phi_i$ such that
\begin{eqnarray}\label{2112}
\lim\phi_i(x)=\phi_i.
\end{eqnarray}
In particular, for an FGM joint distribution, condition (\ref{2111}) reduces to
$$1+\gamma_{12}+\gamma_{13}+\gamma_{23}> 0.$$

For the distribution, under conditions (\ref{2111}) and (\ref{2112}), some direct calculations lead to the following  uniformly asymptotic results over all $s,\ z\ge 0$: for $i=1,2$,
\begin{eqnarray*}\label{221}
P(X_k^{(i)}\in x_i+\Delta^{(i)}|\theta_k=s)
&\sim& F_i(x_i+\Delta^{(i)})\big(1+\gamma_{i3}\phi_3(s)\phi_i\big)=F_i(x_i+\Delta^{(i)})h_i(s),
\end{eqnarray*}
\begin{eqnarray*}\label{224}
P\big(\vec{X}_k\in\vec{x}+\vec{\Delta}|\theta_k=s)
&\sim&\prod_{i=1}^{2} F_i(x_i+\Delta^{(i)})\big(1+\gamma_{12}\phi_1\phi_2+\gamma_{13}\phi_1\phi_3(s)+\gamma_{23}\phi_2\phi_3(s)\big)\nonumber\\
&=&\prod_{i=1}^{2}F_i(x_i+\Delta^{(i)})g(s)
\end{eqnarray*}
and for $1\le i\neq j\le2$,
\begin{eqnarray*}\label{222}
P(X_k^{(i)}\in x_i+\Delta^{(i)}|X_k^{(j)}=z,\theta_k=s)&\sim&
F_i(x_i+\Delta^{(i)})\Big(1+\frac{\gamma_{ij}\phi_i\phi_j(z)+\gamma_{i3}\phi_i\phi_3(s)}
{1+\gamma_{j3}\phi_j(z)\phi_3(s)}\Big)\nonumber\\
&=&F_i(x_i+\Delta^{(i)})g_{ij}(z,s).
\end{eqnarray*}

From Proposition 1.1 of Yang and Wang (2013), there
exists a positive constant $C_i^{\prime}$ such
that $|\phi_i(x)|\leq C_i^{\prime}$ for all $x\in (-\infty,+\infty),\
i=1,2,3$. Hence, for each fixed $T\in\Lambda$,
$$b^*=\max\{\sup_{s\in [0,T]}1+\gamma_{i3}\phi_3(s)\phi_i:\ i=1,2\}<\infty,$$
$$b_*=\min\{\inf_{s\in [0,T]}1+\gamma_{i3}\phi_3(s)\phi_i:\ i=1,2\}>0,$$
$$d^*=\sup_{s\in[0,T)}\big(1+\gamma_{12}\phi_1\phi_2+\gamma_{13}\phi_1\phi_3(s)+\gamma_{23}\phi_2\phi_3(s)\big)<\infty,$$
$$d_*=\inf_{s\in[0,T)}\big(1+\gamma_{12}\phi_1\phi_2+\gamma_{13}\phi_1\phi_3(s)+\gamma_{23}\phi_2\phi_3(s)\big)>0,$$
$$a^*=\max\Big\{\sup_{z\in(-\infty,\infty),\ s\in[0,T)}\Big(1+\frac{\gamma_{ij}\phi_i\phi_j(z)+\gamma_{i3}\phi_i\phi_3(s)}
{1+\gamma_{j3}\phi_j(z)\phi_3(s)}\Big):1\le i\neq j\le2\Big\}<\infty$$
and
$$a_{*}=\min\Big\{\inf_{z\in(-\infty,\infty),\ s\in[0,T)}\Big(1+\frac{\gamma_{ij}\phi_i\phi_j(z)+\gamma_{i3}\phi_i\phi_3(s)}
{1+\gamma_{j3}\phi_j(z)\phi_3(s)}\Big):1\le i\neq j\le2\Big\}>0.$$
\vspace{0.2cm}

\indent{\bf Copula 4.2.} The tri-dimensional Frank copula is of the form
\begin{eqnarray*}
C(u,v,w)=-\frac{1}{\gamma}\ln\Big(1+\frac{(e^{-\gamma u}-1)(e^{-\gamma v}-1)(e^{-\gamma w}-1)}{(e^{-\gamma}-1)^2}\Big),\quad u,v,w\in[0,1],
\end{eqnarray*}
where $\gamma$ is a positive constant.\\

Some direct calculations
lead to the following results: for $i = 1,2$,
\begin{align}\label{001}
&P(X_k^{(i)}\in x_i+\Delta^{(i)}|\theta_k=s)=\lim_{d\downarrow0}\frac{P(X_k^{(i)}\in x_i+\Delta^{(i)},\theta_k\in [s,s+d))}{P(\theta_k\in [s,s+d))}\nonumber\\
&=\frac{e^{-\gamma G(s)}\big(e^{-\gamma}-1)(e^{-\gamma F_i(x_i+d_i)}-e^{-\gamma F_i(x_i)})}{\big((e^{-\gamma}-1)+(e^{-\gamma G(s)}-1)(e^{-\gamma F_i(x_i)}-1)\big)\big((e^{-\gamma}-1)+(e^{-\gamma G(s)}-1)(e^{-\gamma F_i(x_i+d_i)}-1)\big)}
\nonumber\\
&\sim F_i(x_i+\Delta^{(i)})\gamma e^{\gamma G(s)}(e^{\gamma}-1)^{-1}=F_i(x_i+\Delta^{(i)})h_i(s),
\end{align}
\begin{eqnarray*}\label{002}
&&P\Big(\vec{X}_k\in\vec{x}+\vec{\Delta}|\theta_k=s)=\lim_{d\downarrow0}P(\vec{X}_k\in\vec{x}+\vec{\Delta},\theta_k\in [s,s+d))\big(P(\theta_k\in [s,s+d))\big)^{-1}\nonumber\\
&=&\frac{e^{-\gamma F_1(x_1+d_1)}-e^{-\gamma F_1(x_1)}}{(e^{-\gamma}-1)^2+(e^{-\gamma F_1(x_1+d_1)}-1)(e^{-\gamma F_2(x_2+d_2)}-1)(e^{-\gamma G(s)}-1)}
\nonumber\\
& &\cdot\frac{(e^{-\gamma F_2(x_2+d_2)}-e^{-\gamma F_2(x_2)})e^{-\gamma G(s)}}{(e^{-\gamma}-1)^2+(e^{-\gamma F_1(x_1+d_1)}-1)(e^{-\gamma F_2(x_2)}-1)(e^{-\gamma G(s)}-1)}
\nonumber\\
& &\cdot\frac{(e^{-\gamma}-1)^2}{(e^{-\gamma}-1)^2+(e^{-\gamma F_1(x_1)}-1)(e^{-\gamma F_2(x_2+d_2)}-1)(e^{-\gamma G(s)}-1)}\nonumber\\
& &\cdot\Big(\frac{(e^{-\gamma}-1)^4}{(e^{-\gamma}-1)^2+(e^{-\gamma F_1(x_1)}-1)(e^{-\gamma F_2(x_2)}-1)(e^{-\gamma G(s)}-1)}
\nonumber\\
& &-\frac{(e^{-\gamma F_1(x_1+d_1)}-1)(e^{-\gamma F_1(x_1)}-1)(e^{-\gamma F_2(x_2+d_2)}-1)(e^{-\gamma F_2(x_2)}-1)(e^{-\gamma G(s)}-1)^2}
{(e^{-\gamma}-1)^2+(e^{-\gamma F_1(x_1)}-1)(e^{-\gamma F_2(x_2)}-1)(e^{-\gamma G(s)}-1)}\Big)\nonumber\\
&\sim&\prod_{i=1}^{2} F_i(x_i+\Delta^{(i)})\gamma^2(2e^{2\gamma G(s)}-e^{\gamma G(s)})(e^{\gamma}-1)^{-2}=\prod_{i=1}^{2}F_i(x_i+\Delta^{(i)})g(s)
\end{eqnarray*}
and for $1\le i\not=j \le 2$,
\begin{eqnarray*}\label{003}
&&P(X_k^{(i)}\in x_i+\Delta^{(i)}|X_k^{(j)}=z,\theta_k=s)\nonumber\\
&=&\frac{(e^{-\gamma}-1)\big((e^{-\gamma}-1)^4-(e^{-\gamma F_i(x_i)}-1)(e^{-\gamma F_i(x_i+d_i)}-1)(e^{-\gamma F_j(z)}-1)^2(e^{-\gamma G(s)}-1)^2\big)}{\big((e^{-\gamma}-1)^2+(e^{-\gamma F_i(x_i)}-1)(e^{-\gamma F_j(z)}-1)(e^{-\gamma G(s)}-1)\big)^2}
\nonumber\\
& &\cdot\frac{(e^{-\gamma F_i(x_i+d_i)}-e^{-\gamma F_i(x_i)})\big((e^{-\gamma}-1)+(e^{-\gamma F_j(z)}-1)(e^{-\gamma G(s)}-1)\big)^2}
{\big((e^{-\gamma}-1)^2+(e^{-\gamma F_i(x_i+d_i)}-1)(e^{-\gamma F_j(z)}-1)(e^{-\gamma G(s)}-1)\big)^2}
\nonumber\\
&\sim& \frac{\gamma F_i(x_i+\Delta^{(i)})}{e^{\gamma}-1}
\cdot\frac{(e^{-\gamma}-1)-(e^{-\gamma F_j(z)}-1)(e^{-\gamma G(s)}-1)}
{(e^{-\gamma}-1)+(e^{-\gamma F_j(z)}-1)(e^{-\gamma G(s)}-1)}=F_i(x_i+\Delta^{(i)})g_{ij}(z,s),
\end{eqnarray*}
where $g_{ij}(z,s)\le\gamma(2e^{\gamma}-1)(e^{\gamma}-1)^{-1}$.

For each fixed $T\in\Lambda$, by $\gamma>0$ and $G(T)<1$, we have
\begin{equation}\label{004}
0<b_*\le b^*<\infty.
\end{equation}
Further, by $\gamma>0$,
\begin{eqnarray*}\label{005}
&&d^*=d^*(T)=\sup_{s\in[0,T]}\gamma^2(2e^{2\gamma G(s)}-e^{\gamma G(s)})(e^{\gamma}-1)^{-2}<\infty.
\end{eqnarray*}
\begin{eqnarray*}\label{006}
d_*=d_*(T)&=\inf_{s\in[0,T]}\gamma^2(2e^{2\gamma G(s)}-e^{\gamma G(s)})(e^{\gamma}-1)^{-2}>0,
\end{eqnarray*}
\begin{eqnarray*}\label{007}
a^*=\max\Big\{\sup_{z\ge0,s\in[0,T]}\frac{(e^{-\gamma}-1)-(e^{-\gamma F_i(z)}-1)(e^{-\gamma G(s)}-1)}
{(e^{-\gamma}-1)+(e^{-\gamma F_i(z)}-1)(e^{-\gamma G(s)}-1)}:i=1,2\Big\}<\infty
\end{eqnarray*}
and
\begin{eqnarray*}\label{008}
a_*=\min\Big\{\inf_{z\ge0,s\in[0,T]}\frac{(e^{-\gamma}-1)-(e^{-\gamma F_i(z)}-1)(e^{-\gamma G(s)}-1)}
{(e^{-\gamma}-1)+(e^{-\gamma F_i(z)}-1)(e^{-\gamma G(s)}-1)}:i=1,2\Big\}>0.
\end{eqnarray*}
\vspace{0.2cm}

\indent{\bf Copula 4.3.} Recall that the product copula and the
bivariate Frank copula are of the forms,
respectively, $$\Pi(u,v)=uv,\ u,v\in[0,1]$$ and for some
$\gamma>0$,
$$C_{\gamma}(\kappa,w)=-\gamma^{-1}\ln\Big(1+\frac{ (e^{-\gamma \kappa}-1)(e^{-\gamma w}-1)}{e^{-\gamma}-1}\Big),\ \kappa,w\in[0,1].$$
We construct a new tri-dimensional function as follows: for any constant $\gamma>0$,
\begin{eqnarray*}\label{21}
C(u,v,w)=C_{\gamma}(\Pi(u,v),w)=-\frac{1}{\gamma}\ln\Big(1+\frac{ (e^{-\gamma uv}-1)(e^{-\gamma w}-1)}{e^{-\gamma}-1}\Big),\ u,v,w\in[0,1].
\end{eqnarray*}
Now, we show that this function $C(u,v,w)$ is a tri-dimensional copula for $\gamma>0$.

First, it is easy to verify that $C(u,v,w)$ satisfies (2.10.4a) and (2.10.4b)
of Nelsen (2006), that is, $$C(0,v,w)=C(u,0,w)=C(u,v,0)=0$$
and
$$C(u,1,1)=u,C(1,v,1)=v,C(1,1,w)=w.$$
Second, the C-volume of the tri-dimensional function on a rectangle
$B=[u_1,u_2]\times[v_1,v_2]\times[w_1,w_2]$ is given, after some
simplification, by
$$V_C(B)=-\gamma^{-1}\ln\Big(1+\frac{Q_1(u_1,u_2,v_1,v_2,w_1,w_2)}{Q_2(u_1,u_2,v_1,v_2,w_1,w_2)}\Big)
=-\gamma^{-1}\ln\Big(1+\frac{Q_1}
{Q_2}\Big),$$
where
\begin{eqnarray*}
Q_1&=&(e^{-\gamma}-1)^3(e^{-\gamma w_2}-e^{-\gamma w_1})(l_{11}-l_{12}-l_{21}+l_{22})\nonumber\\
& &+(e^{-\gamma}-1)^2(e^{-\gamma w_2}-e^{-\gamma w_1})(e^{-\gamma w_1}+e^{-\gamma w_2}-2)(l_{11}l_{22}-l_{12}l_{21})\nonumber\\
& &+(e^{-\gamma}-1)(e^{-\gamma w_1}-1)(e^{-\gamma w_2}-1)(e^{-\gamma w_2}-e^{-\gamma w_1})l_{11}l_{12}l_{21}l_{22}
(l_{21}^{-1}-l_{22}^{-1}
-l_{11}^{-1}+l_{12}^{-1}),
\end{eqnarray*}
\begin{eqnarray*}
Q_2=\prod_{k=1}^{2}
\big((e^{-\gamma}-1)+l_{kk}(e^{-\gamma w_1}-1)\big)
\prod_{1\le i\not=j \le2}
\big((e^{-\gamma}-1)+l_{ij}(e^{-\gamma w_2}-1)\big)
\end{eqnarray*}
and $l_{ij}=e^{-\gamma u_i v_j}-1$ for $1\le i,j \le2$.
Notice that, by $l_{12}l_{21}>l_{11}l_{22}$ and
$$(1+l_{11})(1+l_{22})<(1+l_{12})(1+l_{21}),$$
\begin{eqnarray*}
Q_1&\leq&(e^{-\gamma}-1)^3(e^{-\gamma w_2}-e^{-\gamma w_1})\big(l_{11}-l_{12}-l_{21}+l_{22}+2l_{11}l_{22}-2l_{12}l_{21}\nonumber\\
& &\ \ \ \ \ \ \ \ \ \ \ \ \ \ \ \ \ \ \ \ \ \ \ \ \ \ \ \ \ \ \ \ \ \
+l_{11}l_{12}l_{21}l_{22}
(l_{21}^{-1}-l_{22}^{-1}
-l_{11}^{-1}+l_{12}^{-1})\big)\nonumber\\
&=&(e^{-\gamma}-1)^3(e^{-\gamma w_2}-e^{-\gamma w_1})
\big((1-l_{12}l_{21})(1+l_{11})(1+l_{22})\nonumber\\
&&\ \ \ \ \ \ \ \ \ \ \ \ \ \ \ \ \ \ \ \ \ \ \ \ \ \ \ \ \ \ \ \ \ \ -(1-l_{11}l_{22})(1+l_{12})(1+l_{21})\big)<0
\end{eqnarray*}
and $Q_2>0$. Hence, $C(u,v,w)$ is a
tri-dimensional copula if and only if $\gamma>0$.

Because $(X^{(i)},\theta)$ has a bivariate Frank copula for $i=1,2$,
by (\ref{001}), we have
\begin{eqnarray*}\label{22}
P(X_k^{(i)}\in x_i+\Delta^{(i)}|\theta_k=s)
&=&\lim_{d\downarrow0}P(X_k^{(i)}\in x_i+\Delta^{(i)},\theta_k\in [s,s+d))\big(P(\theta_k\in [s,s+d)\big)^{-1}\nonumber\\
&\sim& F_i(x_i+\Delta^{(i)})\gamma e^{\gamma G(s)}(e^{\gamma}-1)^{-1}\nonumber\\
&=&F_i(x_i+\Delta^{(i)})h_i(s),
\end{eqnarray*}
\begin{eqnarray*}\label{24}
&&P(\vec{X}_k\in\vec{x}+\vec{\Delta}|~\theta_k=s)=\lim_{d\downarrow0}P(\vec{X}_k\in\vec{x}+\vec{\Delta},\theta_k\in [s,s+d))\big(P(\theta_k\in [s,s+d))\big)^{-1}\nonumber\\
&=&(e^{-\gamma}-1)\big((e^{-\gamma}-1)+(e^{-\gamma F_1(x_1)F_2(x_2)}-1)(e^{-\gamma G(s)}-1)\big)^{-1}\nonumber\\
&&\ \ \ \ \ \ \ \ \ \ \ \ \cdot\big((e^{-\gamma}-1)+(e^{-\gamma F_1(x_1+d_1)F_2(x_2+d_2)}-1)(e^{-\gamma G(s)}-1)\big)^{-1}\nonumber\\
& &\cdot e^{-\gamma G(s)}\big((e^{-\gamma}-1)+(e^{-\gamma F_1(x_1)F_2(x_2+d_2)}-1)(e^{-\gamma G(s)}-1)\big)^{-1}\nonumber\\
&&\ \ \ \ \ \ \ \ \ \ \ \ \cdot\big((e^{-\gamma}-1)+(e^{-\gamma F_1(x_1+d_1)F_2(x_2)}-1)(e^{-\gamma G(s)}-1)\big)^{-1}\nonumber\\
& &\cdot\Big((e^{-\gamma}-e^{-\gamma G(s)})^2(e^{-\gamma F_1(x_1)F_2(x_2)}-e^{-\gamma F_1(x_1)F_2(x_2+d_2)}-e^{-\gamma F_1(x_1+d_2)F_2(x_2)}\nonumber\\
&&\ \ \ \ \ \ \ \ \ \ \ \ +e^{-\gamma F_1(x_1+d_1)F_2(x_2+d_2)})\nonumber\\
& &\ \ \ +2(e^{-\gamma}-e^{-\gamma G(s)})(e^{-\gamma G(s)}-1)e^{-\gamma (F_1(x_1)F_2(x_2)+F_1(x_1+d_1)F_2(x_2+d_2))}\nonumber\\
&&\ \ \ \ \ \ \ \ \ \ \ \ \cdot(1-e^{\gamma \prod_{i=1}^{2}F_i(x_i+\Delta^{(i)})})\nonumber\\
& &\ \ \ +(e^{-\gamma G(s)}-1)^2e^{-\gamma\big(F_1(x_1)F_2(x_2)+F_1(x_1+d_1)F_2(x_2+d_2)
+F_1(x_1)F_2(x_2+d_2)+F_1(x_1+d_1)F_2(x_2)\big)}\nonumber\\
& &\ \ \ \ \ \ \ \ \ \ \ \ \cdot(e^{\gamma F_1(x_1+d_1)F_2(x_2)}-e^{\gamma F_1(x_1)F_2(x_2)}-e^{\gamma F_1(x_1+d_1)F_2(x_2+d_2)}+e^{\gamma F_1(x_1)F_2(x_2+d_2)})\Big)\nonumber\\
&\sim&\prod_{i=1}^{2}F_i(x_i+\Delta^{(i)})\big(\gamma (e^{-\gamma}-e^{-2\gamma})e^{\gamma G(s)}+\gamma^2 (e^{-\gamma}+e^{-2\gamma})(e^{\gamma G(s)}-1)\big)(1-e^{-\gamma})^{-2}\nonumber\\
&=&\prod_{i=1}^{2}F_i(x_i+\Delta^{(i)})g(s),
\end{eqnarray*}
and for $1\le i\not=j \le 2$,
\begin{eqnarray*}\label{241}
&&P(X^{(i)}_k\in x_i+\Delta^{(i)}|X^{(j)}_k=z,\theta_k=s_k)=\big((e^{-\gamma}-1)+(e^{-\gamma F_j(z)}-1)(e^{-\gamma G(s)}-1)\big)^2\nonumber\\
&&\ \ \ \ \ \ \ \ \ \ \ \ \ \ \cdot\big((e^{-\gamma}-1)+(e^{-\gamma F_i(x_i)F_j(z)}-1)(e^{-\gamma G(s)}-1)\big)^{-2}\nonumber\\
&&\ \ \ \ \ \ \ \ \ \ \ \ \ \ \cdot\big((e^{-\gamma}-1)+(e^{-\gamma F_i(x_i+d_i)F_j(z)}-1)(e^{-\gamma G(s)}-1)\big)^{-2}\nonumber\\
& &\ \ \ \ \ \ \ \ \ \ \ \ \ \ \cdot \Big(F_i(x_i+d_i)e^{\gamma \overline {F_i}(x_i+d_i)F_j(z)}(e^{-\gamma}-e^{-\gamma G(s)}+(e^{-\gamma G(s)}-1)e^{-\gamma F_i(x_i)F_j(z)})^2\nonumber\\
& &\ \ \ \ \ \ \ \ \ \ \ \ \ \ \ \ \ \ \ \ \ \ \ \ -F_i(x_i)e^{\gamma \overline {F_i}(x_i)F_j(z)}(e^{-\gamma}-e^{-\gamma G(s)}+(e^{-\gamma G(s)}-1)e^{-\gamma F_i(x_i+d_i)F_j(z)})^2\Big)\nonumber\\
&\sim& F_i(x_i+\Delta^{(i)})\Big(1+\frac{\gamma F_j(z)((1-e^{-\gamma})+(e^{-\gamma F_j(z)}+1)(e^{-\gamma G(s)}-1))}{(e^{-\gamma}-1)+(e^{-\gamma F_j(z)}-1)(e^{-\gamma G(s)}-1)}\Big)\nonumber\\
&=&F_i(x_i+\Delta^{(i)})g_{ij}(z,s),
\end{eqnarray*}
where $g(s)\leq\big(\gamma+\gamma^2(1+e^{-\gamma})\big)(1-e^{-\gamma})^{-1}$ and $g_{ij}(z,s)\leq 1+\gamma(1+e^{\gamma})$.

Similar to copula 4.2 for each fixed $T\in\Lambda$,
due to the fact $\gamma>0$ and $G(T)<1$, we know that
(\ref{004}) holds too. In addition, when $\gamma>0$,
\begin{eqnarray*}\label{2403}
d^*=d^*(T)=\sup_{s\in[0,T]}\frac{\gamma (e^{-\gamma}-e^{-2\gamma})e^{\gamma G(s)}+\gamma^2 (e^{-\gamma}+e^{-2\gamma})(e^{\gamma G(s)}-1)}{(1-e^{-\gamma})^2}<\infty,
\end{eqnarray*}
\begin{eqnarray*}\label{2404}
d_*=d_*(T)&=\inf_{s\in[0,T]}\frac{\gamma (e^{-\gamma}-e^{-2\gamma})e^{\gamma G(s)}+\gamma^2 (e^{-\gamma}+e^{-2\gamma})(e^{\gamma G(s)}-1)}{(1-e^{-\gamma})^2}>0,
\end{eqnarray*}
\begin{eqnarray*}\label{2401}
a^*=\sup_{z\ge0,s\in[0,T]}\Big(1+\frac{\gamma F_j(z)\big((1-e^{-\gamma})+(e^{-\gamma F_j(z)}+1)(e^{-\gamma G(s)}-1)\big)}{(e^{-\gamma}-1)+(e^{-\gamma F_j(z)}-1)(e^{-\gamma G(s)}-1)}\Big)<\infty,
\end{eqnarray*}
and when $0<\gamma<1$,
\begin{eqnarray*}\label{2402}
a_*=\inf_{z\ge0,s\in[0,T]}\Big(1+\frac{\gamma F_j(z)\big((1-e^{-\gamma})+(e^{-\gamma F_j(z)}+1)(e^{-\gamma G(s)}-1)\big)}{(e^{-\gamma}-1)+(e^{-\gamma F_j(z)}-1)(e^{-\gamma G(s)}-1)}\Big)>0.
\end{eqnarray*}

\setcounter{equation}{0}
\setcounter{lemma}{0}
\section{A example}

In this section, we give a non locally almost decreased distribution in the class $\mathcal{S}_{loc}$.

Let $\{a_0=0,a_n=2^{{n^2}},b_n=a_n+a_{m_n}\ln^2(n+1):n\ge1\}$ be a sequence of positive numbers, where $m_n=\min\{k:k\ge\sqrt{5}(\sqrt{6})^{-1}n\}$ for all $n\ge1$. Clearly,
\begin{eqnarray}\label{500}
\sqrt{5}(\sqrt{6})^{-1}n\le m_n<\sqrt{5}(\sqrt{6})^{-1}(n+1)+1.
\end{eqnarray}
And let $f_0:[0,\infty)\longmapsto (0,1]$ be a linear function such that
\begin{eqnarray}\label{501}
&&f_0(x)=\big(f_0(a_0)+\frac{f_0(a_1)-f_0(a_0)}{a_1-a_0}(x-a_0)\big)\textbf{1}(a_0\le x\le a_1)\nonumber\\
&&\ \ \ \ \ \ \ \ \ \ \ \ \ \ +
\sum_{n=1}^\infty\Big(\big(f_0(a_n)+\frac{f_0(b_n)-f_0(a_n)}{b_n-a_n}(x-a_n)\big)\textbf{1}(a_n<x\le b_n)\nonumber\\
&&\ \ \ \ \ \ \ \ \ \ \ \ \ \ +\big(f_0(b_n)+\frac{f_0(\frac{a_{n+1}+b_n}{2})-f_0(b_n)}{\frac{a_{n+1}+b_n}{2}-b_n}(x-b_n)\big)\textbf{1}(b_n<x\le \frac{a_{n+1}+b_n}{2})\nonumber\\
&&\ \ \ \ \ \ \ \ \ \ \ \ \ \ +\big(f_0(\frac{a_{n+1}+b_n}{2})+\frac{f_0(a_{n+1})-f_0(\frac{a_{n+1}+b_n}{2})}{\frac{a_{n+1}+b_n}{2}-b_n}(x-\frac{a_{n+1}+b_n}{2})\big)\nonumber\\
&&\ \ \ \ \ \ \ \ \ \ \ \ \ \ \ \ \ \ \ \ \ \ \ \ \ \ \ \ \cdot\textbf{1}(\frac{a_{n+1}+b_n}{2}<x\le a_{n+1})\Big),
\end{eqnarray}
$f_0(a_n)=f_0(\frac{a_{n+1}+b_n}{2})=2a_n^{-3}$ and $f_0(b_n)=\frac{f_0(a_n)}{\ln(n+1)}$, for all $n\ge1$.

Because
$$\int_{a_n}^{a_{n+1}} f_0(y)dy\le f_0(a_n)a^{2+\alpha}_{n+1}=2^{-n^2+2(2n+1)}$$
for all $n\ge1$, $0<a=\int_{0}^{\infty} f_0(x)dx<\infty$. Therefore, the function $f=a^{-1}f_0$ is a density corresponding to a distribution $F$ supported on $[0,\infty)$. Without loss of generality, set $a=1$.

The density $f$ is not almost decreased. In fact, when $n\to\infty$,
$$f(b_n)=f(a_n)\big(\ln(n+1)^{-1}\big)=o\big(f(a_n)\big)=o\big(f(\frac{a_{n+1}+b_n}{2})\big).$$
And then by (\ref{501}), when $n\to\infty$, we have
$$\frac{f(a_n)-f(b_n)}{b_n-a_n}=\frac{f(a_n)\big(1-\ln^{-1}(n+1)\big)}{a_{m_n}\ln^2(n+1)}=o\big(\frac{f(b_n)}{a_{m_n}}\big),$$
$$\frac{f(\frac{a_{n+1}+b_n}{2})-f(b_n)}{\frac{a_{n+1}+b_n}{2}-b_n}\sim\frac{2f(a_n)\big(1-\ln^{-1}(n+1)\big)}{a_{n+1}}
=o\big(\frac{f(b_n)\ln^2(n+1)}{a_{n+1}}\big),$$
and
$$\frac{f(\frac{a_{n+1}+b_n}{2})-f(a_{n+1})}{a_{n+1}-\frac{a_{n+1}+b_n}{2}}\sim\frac{2f(a_n)-f(a_{n+1})}{a_{n+1}}=o\big(f(a_{n+1})a^{-1}_nn^{5n}\big).$$
Thus,
\begin{eqnarray}\label{502}
\sup_{y\in J_{ni}}\mid f^{\prime}(y)\mid=o\big(\inf_{x\in J_{ni}}f(x)\big),\ i=1,2,3,
\end{eqnarray}
where $J_{n1}=(a_n,b_n],\ J_{n2}=(b_n,\frac{a_{n+1}+b_n}{2}]$, and $J_{n3}=(\frac{a_{n+1}+b_n}{2},a_{n+1}]$.

For $n\ge1$, we denote any two adjacent numbers in set $\{a_n,b_n,\frac{a_{n+1}+b_n}{2},a_{n+1}\}$ by $c_n,d_n$. By the method of Lemma 4.1 in Xu et al. (2015) and (\ref{502}), we know that, for any fixed number $t>0$ and any $x>0$, there is a $n\ge1$ such that
\begin{eqnarray*}
f(x)&=&f(x-t)+\int_{x-t}^xf^{\prime}(y)dy\textbf{1}(c_n<x-t<x\le d_n)\nonumber\\
&&\ \ \ \ \ \ \ \ \ \ \ \ \ \ +\Big(\int_{x-t}^{d_n}+\int_{d_n}^x\Big)f^{\prime}(y)dy\textbf{1}(c_n<x-t\le d_n<x)\nonumber\\
&\sim& f(x-t).
\end{eqnarray*}
Thus, the density $f$ belongs to the long-tailed function class, denoted by $f\in\mathcal{L}_d$. Further, the corresponding distribution $F\in\mathcal{L}_{loc}$, see Asmussen et al. (2003).


Now, for all $n\ge1$ and $x\in J_n=(a_n,a_{n+1}]$, when $n\to\infty$, or equivalently $x\to\infty$, we are going to prove that
\begin{eqnarray}\label{507}
I(x)&=&\int_0^x f(x-y)f(y)dy=\Big(2\int_0^{a_{m_n}}+ \int_{a_{m_n}}^{x-a_{m_n}}\Big) f(x-y)f(y)dy\nonumber\\
&=&2I_1(x)+I_2(x)\sim 2f(x),
\end{eqnarray}
that is the density $f$ belongs to the subexponential function class, denoted by $f\in\mathcal{S}_d$. Thus, the corresponding distribution $F\in\mathcal{S}_{loc}$, see also Asmussen et al. (2003).

For $ x\in(a_n,a_{n+1}]$, because $y$ and $x-y\in(a_{m_n},x-a_{m_n}]$, by (\ref{501}), we have
\begin{eqnarray}\label{508}
I_2(x)&\le& f^2(a_{m_n})a_{n+1}=2^{-6m_n^2+(n+1)^2+2}
=2^{-3(n+1)^2+3(n+1)^2-5n^2+(n+1)^2+2}\nonumber\\
&=&o\big(f(a_{n+1})\big)=o\big(f(x)\big).
\end{eqnarray}

In the following, we deal with $I_1(x)$. Because $f\in\mathcal{L}_d$, we just have to prove that
\begin{eqnarray}\label{5080}
f(x-a_{m_n})\sim f(x)
\end{eqnarray}
for $x\in J_{ni},i=1,2,3$, respectively.

When $x\in J_{n1}=(a_n,b_n]$, because $0\le y\le a_{m_n}$,
$$\frac{a_{n}+b_{n-1}}{2}<a_n-a_{m_n}<x-a_{m_n}\le x-y\le x\le b_n.$$
If $a_n\le x-a_{m_n}<x\le b_n$, then by (\ref{501}) and (\ref{500}),
\begin{eqnarray*}
f(x-a_{m_n})&=&f(x)+\frac{f(a_n)\big(1-\ln^{-1}(n+1)\big)}{a_{m_n}\ln^2(n+1)}\nonumber\\
&=&f(x)+\frac{f(b_n)\big(1-\ln^{-1}(n+1)\big)}{a_{m_n}\ln(n+1)}\sim f(x).
\end{eqnarray*}
If $\frac{a_{n}+b_{n-1}}{2}< a_n-a_{m_n}\le x-a_{m_n}<a_n<x\le b_n$, then by (\ref{501}) and (\ref{500}),
\begin{eqnarray*}
f(a_n)&\le& f(x-a_{m_n})\le f(a_n-a_{m_n})\nonumber\\
&=&f(a_{n})+\frac{2\big(f(a_{n-1})-f(a_n)\big)a_{m_n}}{a_n-b_{n-1}}\nonumber\\
&=&f(a_{n})+f(a_{n})2^{-\frac{1}{6}n^2+6n}\sim f(a_n),
\end{eqnarray*}
that is $f(x-a_{m_n})\sim f(a_n)$. On other hand, because $x-a_n\le a_{m_n}$, we have
\begin{eqnarray*}
f(x)=f(a_n)+\frac{f(a_n)\big(\ln^{-1}(n+1)-1\big)}{a_{m_n}\ln^2(n+1)}\cdot (x-a_n)\sim f(a_n),
\end{eqnarray*}
that is (\ref{5080}) holds for $x\in J_{n1}=(a_n,b_n]$.

When $x\in J_{n2}=(b_n,\frac{a_{n+1}+b_{n}}{2}]$, then
$$a_{n}<b_n-a_{m_n}<x-a_{m_n}<x\le \frac{a_{n+1}+b_{n}}{2}.$$
If $b_{n}<x-a_{m_n}<x\le \frac{a_{n+1}+b_{n}}{2},$ then by (\ref{501}) and (\ref{500}), we have
$$f(x-a_{m_n})=f(x)-\frac{2\big(f(a_n)-f(b_n)\big)}{a_{n+1}-b_n}\cdot a_{m_n}\sim f(x).$$
If $a_n<b_n-a_{m_n}<x-a_{m_n}\le b_n<x,$ then by (\ref{501}), we have
$$f(x-a_{m_n})\le f(b_n-a_{m_n})=f(b_n)+\frac{f(a_n)\big(1-\ln^{-1}(n+1)\big)}{a_{m_n}\ln^2(n+1)}\cdot a_{m_n}\sim f(b_n).$$
On other hand, because $x-b_n\le a_{m_n}$, we have
\begin{eqnarray*}
f(x)=f(b_n)+\frac{2\big(f(a_n)-f(b_n)\big)}{a_{n+1}-b_n}\cdot (x-b_n)\sim f(b_n)\le f(x-a_{m_n}).
\end{eqnarray*}
Thus, (\ref{5080}) holds for $x\in(b_n,\frac{a_{n+1}+b_{n}}{2}]$.

When $x\in J_{n3}=(\frac{a_{n+1}+b_{n}}{2},a_{n+1}]$, then
$$\frac{a_{n+1}+b_{n}}{2}-a_{m_n}<x-a_{m_n}<x\le a_{n+1}.$$
If $\frac{a_{n+1}+b_{n}}{2}<x-a_{m_n}<x\le a_{n+1}$, then
$$f(x-a_{m_n})=f(x)+\frac{2\big(f(a_n)-f(a_{n+1})\big)}{a_{n+1}-b_n}\cdot a_{m_n}\sim f(x).$$
If $b_n<\frac{a_{n+1}+b_{n}}{2}-a_{m_n}<x-a_{m_n}\le\frac{a_{n+1}+b_{n}}{2}<x\le a_{n+1}$, then by (\ref{501}),
$$f(x-a_{m_n})\ge f(\frac{a_{n+1}+b_{n}}{2}-a_{m_n})=f(a_n)+\frac{2\big(f(a_n)-f(b_n)\big)}{a_{n+1}-b_n}\cdot a_{m_n}\sim f(a_n);$$
and by $x-\frac{a_{n+1}+b_{n}}{2}\le a_{m_n}$,
$$f(x)=f(a_n)+\frac{2\big(f(a_n)-f(a_{n+1})\big)}{a_{n+1}-b_n}\cdot \big(x-\frac{a_{n+1}+b_{n}}{2}\big)\sim f(a_n)\ge f(x-a_{m_n}).$$
Thus, (\ref{5080}) holds for $x\in(\frac{a_{n+1}+b_{n}}{2},a_{n+1}]$.
\\
\\

\textbf{Acknowledgements} The authors are very grateful to Professor Dmitry Korshunov and Dr. Changjun Yu for their helpful discussions and comments.

\end{document}